\documentclass{article}[10pt,oneside,a4paper,draft]

\usepackage{amsmath}
\usepackage{amssymb}
\usepackage{tikz}
\usepackage[footnotesize,hang]{caption}
\usepackage{pstricks}
\usepackage{subfig}
\usepackage{color}
\usepackage{rotating}
\usepackage{mathtools}
\usepackage{xcolor}
\usepackage{pifont}
\usepackage{graphicx}
\usepackage{authblk}

\usetikzlibrary{decorations.pathmorphing}
\renewcommand{\inf}{\infty}
\newcommand{\eps}{\epsilon}

\newcommand{\e}{\mathrm{e}}

\renewcommand{\d}{\,\mathrm{d}}
\newcommand{\diff}[2]{\frac{\mathrm{d} #1}{\mathrm{d} #2}}

\newcommand{\bif}{\lambda}
\newcommand{\xoneeqm}{\tilde{x}}
\newcommand{\xtwoeqm}{\tilde{\tilde{x}}}
\newcommand{\ie}{\emph{i.e.}}
\newcommand{\eg}{\emph{e.g.}}
\newcommand{\Ord}{\mathcal{O}}
\newcommand{\ord}{\operatorname{ord}}
\newcommand{\erfi}{\operatorname{erfi}}
\newcommand{\erf}{\operatorname{erf}}
\newcommand{\ODeltaE}{O$\mathrm{\Delta}$E}
\newcommand{\PDeltaE}{P$\mathrm{\Delta}$E}
\newcommand{\DeltaDE}{$\mathrm{\Delta}$DE}

\def\XXint#1#2#3{{\setbox0=\hbox{$#1{#2#3}{\int}$}
\vcenter{\hbox{$#2#3$}}\kern-.5\wd0}}

\usepackage{setspace}
\begin{document}

\title{Multiple scales and matched asymptotic expansions for the discrete logistic equation}

\author[ ]{Cameron L.~Hall\thanks{Electronic address: \texttt{hall@maths.ox.ac.uk}; Corresponding author}}
\author[ ]{Christopher J.~Lustri\thanks{Electronic address: \texttt{christopher.lustri@sydney.edu.au}}}
\affil[*]{Mathematical Institute, University of Oxford, Woodstock Rd, Oxford, {OX2 6GG}, United Kingdom. hall@maths.ox.ac.uk.}
\affil[$\dag$]{Department of Mathematics and Statistics, Faculty of Science, University of Sydney, Sydney, 2006, Australia. christopher.lustri@sydney.edu.au.}

%\author{Cameron L.~Hall         \and
%        Christopher J.~Lustri

%
%\institute{C. L. Hall \at
              %Mathematical Institute\\
              %University of Oxford\\
              %Woodstock Rd\\
              %Oxford OX2 6GG\\
              %United Kingdom\\
              %Tel.: +44-1865-280618\\
              %\email{hall@maths.ox.ac.uk}
           %\and
           %C. J. Lustri \at
           %Faculty of Science\\
           %Carslaw Building\\
           %University of Sydney\\
           %New South Wales 2006\\
           %Australia\\
           %Tel.: +61-2-9351-3879\\
           %\email{christopher.lustri@sydney.edu.au}
%% The correct dates will be entered by the editor
%
\maketitle

\begin{abstract}
In this paper, we combine the method of multiple scales and the method of matched asymptotic expansions to construct uniformly-valid asymptotic solutions to autonomous and non-autonomous difference equations in the neighbourhood of a period-doubling bifurcation. In each case, we begin by constructing multiple scales approximations in which the slow time scale is treated as a continuum variable, leading to difference-differential equations. The resultant approximations fail to be asymptotic at late time, due to behaviour on the slow time scale, it is necessary to eliminate the effects of the fast time scale in order to find the late time rescaling, but there are then no difficulties with applying the method of matched asymptotic expansions. The methods that we develop lead to a general strategy for obtaining asymptotic solutions to singularly-perturbed difference equations, and we discuss clear indicators of when multiple scales, matched asymptotic expansions, or a combined approach might be appropriate.   
%\keywords{difference equations \and asymptotic analysis \and discrete-to-continuum asymptotics \and logistic difference equation}
%\subclass{39A99}
\end{abstract}

\allowdisplaybreaks[0]

\section{Introduction}
\label{S:Introduction}

\subsection{Multiple scales for difference equations}

The method of multiple scales is extremely well established as an asymptotic method for analysing differential equations. As early as 1973, Nayfeh \cite{NayfehPert} was able to describe how the method of multiple scales had been  applied to nonlinear oscillator theory, orbital mechanics, flight mechanics, buckling analysis, flutter, wave propagation (in both solids and fluids), dispersive waves, plasma physics, atmospheric science, and statistical mechanics (amongst other areas). The method of multiple scales has continued to be used extensively in the analysis of differential equations, and is covered in most standard texts and courses on perturbation methods (see, for example, \cite{BenderOrszag,HinchPert,HolmesPert,KevorkianCole,NayfehPert}).

In contrast, the application of the method of multiple scales to difference equations is much less well-developed. While Hoppensteadt and Miranker \cite{Hoppensteadt1977} introduced the method of multiple scales for difference equations as early as 1977, there were relatively few developments in this area until the rediscovery of multiple scales (and the Poincar{\'e}--Lindstedt method) for difference equations by Luongo \cite{Luongo1996} and Maccari \cite{Maccari1999} in the 1990s, followed by the innovative recent work by van Horssen and ter Brake \cite{VanHorssen2009}, and Rafei and van Horssen \cite{Rafei2010,Rafei2012,Rafei2014}. Very few texts on perturbation methods discuss difference equations in detail (\cite{BenderOrszag} and \cite{HolmesPert} being two notable exceptions), and Holmes \cite{HolmesPert} recently observed that there are still many open problems on the application of multiple scales to difference equations.

The essence of the method of multiple scales lies in replacing the original independent variable (or system of variables) with a new system that reflects the different scales over which important phenomena occur. In the analysis of the ordinary differential equation (ODE) for a weakly nonlinear oscillator, for example, the method of multiple scales involves transforming the ODE into a partial differential equation (PDE) with two independent time variables: one representing the fast time scale on which the oscillations take place, and the other representing the slow time scale on which the phase and amplitude of the oscillations evolve. %The extra `freedom' associated with the presence of two time variables is then used to eliminate the loss of asymptoticity associated with the secular growth of correction terms over the slow time scale.

The method of multiple scales for an ordinary difference equation (\ODeltaE) is essentially the same in that the original, discrete independent variable is replaced by multiple independent variables that reflect the different time scales inherent to the problem. However, there are two different approaches to slow time variables for difference equations in the present literature. One approach, introduced by Hoppensteadt and Miranker \cite{Hoppensteadt1977}, 
%and also used elsewhere (see, for example, \cite{Luongo1996,Maccari1999,Marathe2006}) 
involves treating the slow time variable as continuous and using Taylor expansions in the continuum variable to transform the original \ODeltaE{} into a difference-differential equation (\DeltaDE). The other approach, introduced by Subramanian and Krishnan \cite{Subramanian1979} (but corrected by van Horssen and ter Brake \cite{VanHorssen2009}), involves treating the slow time variable as discrete, so that the \ODeltaE{} is transformed into a partial difference equation (\PDeltaE). Both approaches continue to be in current use: for applications involving a continuum slow time variable, see for example \cite{Joshi2015,Luongo1996,Maccari1999,Marathe2006,Mickens1987,Mickens2014}; for applications involving a discrete slow time variable, see for example \cite{Rafei2010,Rafei2012,Rafei2014}. %(Note that while Mickens \cite{Mickens1987,Mickens2014} describes the slow time variable as `discrete' throughout, the use of slow time \emph{differential} operators rather than slow time \emph{difference} operators means that the method in these papers is more akin to the Hoppensteadt--Miranker approach than to the van Horssen--ter Brake approach.)

For a weakly nonlinear multiple scales problem, the leading-order problem in the slow time variable will often be nonlinear, even if the leading-order problem in the fast time variable is linear. Since a wider range of tools exist for solving nonlinear differential equations than nonlinear difference equations, it would seem that there are clear advantages to the Hoppensteadt--Miranker approach of treating the slow time variable as continuous. However, difference equations -- especially nonlinear difference equations -- can behave very differently from their analogous differential equations (see \cite{VanHorssen2009} for examples), and important features of a difference equation might be lost if a continuum slow time variable is introduced. This motivates van Horssen and ter Brake \cite{VanHorssen2009} to avoid continuum variables altogether, instead introducing a discrete slow time scale so that the discrete character of the problem can be maintained throughout.

In this paper, we demonstrate by examples that a combination of the method of multiple scales with the method of matched asymptotic expansions can capture the subtleties of a nonlinear \ODeltaE{} without the need for a discrete slow time variable and the challenges that this introduces. In \cite{VanHorssen2009}, van Horssen and ter Brake are correct to observe that the \ODeltaE{} obtained from discretising an ODE may behave very differently from the original ODE (\eg{} the ODE may exhibit finite-time blow-up while the \ODeltaE{} has a solution for all time). However, such discrepancies must also be associated with an asymptotic failure of the original multiple scales expansion, suggesting that they could be corrected using the method of matched asymptotic expansions. (For further discussion of this point, see Section \ref{S:Strategy}.)

We illustrate our method of combining the method of multiple scales with the method of matched asymptotic expansions by investigating the discrete logistic equation in the neighbourhood of a period doubling bifurcation, although the methods that we develop are more broadly applicable to singularly perturbed difference equations. %By concentrating our analysis on a period doubling bifurcation, we can explore in detail how multiple scales approximations can fail, and how matched asymptotic expansions can be used to obtain uniformly valid asymptotic expansions that go beyond the point where the scalings used to obtain the original multiple scales approximation are no longer valid. The methods described in this paper could easily be applied to other singular perturbation problems in difference equations. %especially those associated with the behaviour of a difference equation in the vicinity of a period doubling bifurcation.

\subsection{The discrete logistic equation}
\label{S:IntroDiscLog}

The discrete logistic equation takes the form
\begin{equation}
 x(n+1) = \bif \, x(n) \, \big[1 - x(n)\big], \qquad 0 < x(0) < 1,
 \label{DiscreteLogistic}
\end{equation}
where $\bif$ is a dimensionless parameter with $0 < \bif \leq 4$. The discrete logistic equation is well-known as a model of population dynamics with discrete generations, and it is commonly used as an archetypal example of a nonlinear difference equation that exhibits period-doubling to chaos as the bifurcation parameter, $\bif$, is increased towards 4. %In the present work, we will only be concerned with the first period-doubling bifurcation, long before the onset of chaotic behaviour. 
For a detailed discussion of the bifurcation structure of the discrete logistic equation, see, for example, \cite{May1976,MurrayMathBiol1,StrogatzNonlinear}.

\begin{figure}
\centering
\includegraphics[draft=false,width=0.9\columnwidth]{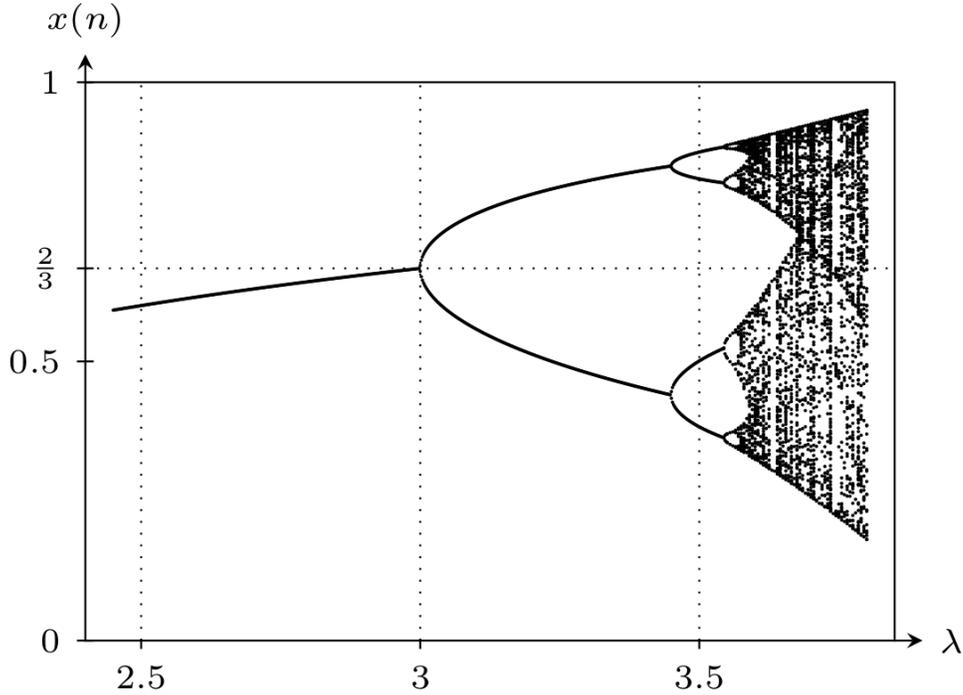}
\caption{Bifurcation diagram showing the stable equilibria of the discrete logistic equation \eqref{DiscreteLogistic} when the bifurcation parameter is varied between $\lambda = 2.45$ and $\lambda = 3.8$.}\label{plot_Bif}
\end{figure}

%In the present work, we focus our attention on the first period-doubling bifurcation of \eqref{DiscreteLogistic}, which occurs at $\bif = 3$. 
For all $\bif > 1$, the discrete logistic equation has an equilibrium solution
\begin{equation}
 \xoneeqm(\bif) =  \frac{\bif-1}{\bif}.
 \label{1PeriodicSol}
\end{equation} 
When $1 < \bif < 3$, this equilibrium is stable. However, as $\bif$ increases through 3, this 1-periodic equilibrium becomes unstable in favour of a 2-periodic cycle that alternates between the values
\begin{equation}
 \xtwoeqm^{\pm}(\bif)= \frac{\bif + 1 \pm \sqrt{\bif^2 - 2 \bif -3}}{2 \bif}.
 \label{2PeriodicSol}
\end{equation} 
For $3 < \bif < 1 + \sqrt{6}$, this 2-periodic cycle is stable, and all trajectories with $0 < x(0) < 1$ and $x(0) \neq \xoneeqm(\bif)$ converge to a cycle which alternates between $\xtwoeqm^+$ and $\xtwoeqm^-$ as $n \rightarrow \infty$. At $\bif = 1 + \sqrt{6}$, the 2-periodic cycle becomes unstable as a stable 4-periodic cycle appears, and the period doubling process continues as shown in Figure \ref{plot_Bif}.

In general, it is difficult to analyse solution trajectories of the discrete logistic equation using asymptotic methods since \eqref{DiscreteLogistic} will be strongly nonlinear for most choices of $\bif$ and $x(0)$. However, by picking $x(0)$ to be close to $\xoneeqm(\bif)$ (or indeed, any periodic solution), \eqref{DiscreteLogistic} can be rescaled as a weakly nonlinear difference equation that is amenable to the method of multiple scales. %If we take $x(0) \approx \xoneeqm(\bif)$ and additionally choose $\bif = 3 + \eps$ for some $0 < \eps \ll 1$, we find that the solution trajectory will initially move away from the unstable $\xoneeqm(\bif)$, until it is captured by the stable cycle represented by $\xtwoeqm^{\pm}(\bif)$. 

In Section \ref{S:StaticAlgebraic}, we take $\bif = 3 + \eps$, and select $x(0)$ near the unstable equilibrium \eqref{1PeriodicSol}. We find that the `plain' method of multiple scale only sees the escape from the unstable equilibrium at $\xoneeqm$, and that matched asymptotic expansions are required in order to describe the approach to the stable cycle at $\xtwoeqm^{\pm}$. Then, in Section \ref{S:Dynamic} we broaden our analysis to consider a non-autonomous nonlinear difference equation: the discrete logistic equation with a slowly varying bifurcation parameter. %This introduces a number of complexities, but we find that the methods introduced in the earlier sections of the paper can again be used to obtain a uniformly valid asymptotic solution.

In \cite{Baesens1991}, Baesens uses the method of renormalisation to perform a comprehensive analysis of the discrete logistic equation with a slowly varying bifurcation parameter. Baesens considers a number of complications that are beyond the scope of the present work: \cite{Baesens1991} deals with backward sweep as well as forward sweep, gives results for the number of observable doublings through the entire period-doubling cascade, and contains a detailed analysis of the critical effects of noise. Although we will discuss the results in \cite{Baesens1991} and draw on some important ideas (such as the stable and unstable adiabatic manifolds), our approach stands independently, focusing on how the asymptotic methods of multiple scales and matched asymptotic expansions can be applied to non-autonomous difference equations. %stands independently of the work in \cite{Baesens1991}, and our focus is on demonstrating that a combination of multiple scales and matched asymptotic expansions can deal with non-autonomous difference equations. %While \cite{Baesens1991} develops a detailed theory for the discrete logistic equation with a slowly varying bifurcation parameter based on renormalisation, our analysis demonstrates that a combined approach of multiple scales and matched asymptotic expansions can also be applied to this problem, illustrating the versatility of this hybrid approach. %and hence our methods are not limited to the analysis of autonomous difference equations.

%However, our analysis in section \ref{S:Dynamic} shows that key features of the discrete logistic equation with a slowly varying bifurcation parameter are also accessible using multiple scales and matched asymptotics, and that these techniques are not limited to the analysis of autonomous difference equations.

% We wish to use the method of multiple scales and the method of matched asymptotic expansions to analyse solution trajectories for the discrete logistic equation when $\bif$ is near the critical value for a period doubling bifurcation.  
% However, by picking $\bif$ so that the system is close to a bifurcation and simultaneously choosing initial conditions near the unstable equilibrium, we can rescale \eqref{DiscreteLogistic} to obtain a weakly nonlinear (but singularly perturbed) difference equation.

\subsection{Outline of paper}
\label{S:Outline}

In this paper, we proceed by considering two examples of the discrete logistic equation where multiple scales and matched asymptotic expansions can be used. These illustrative examples will be of increasing complexity, and they show clearly how to deal with many of the challenges associated with applying multiple scales and matched asymptoic expansions to singularly perturbed difference equations.

In Section \ref{S:StaticAlgebraic}, we demonstrate our method of combining multiple scales with matched asymptotic expansions by analysing the discrete logistic equation with $\bif = 3 + \eps$ and $x(0) = \frac{2}{3}$, where $0 < \eps \ll 1$. Rescaling $x(n)$ in \eqref{DiscreteLogistic} by introducing $x(n) = \frac{2}{3} + \eps\, X(n)$, we obtain a weakly nonlinear difference equation that gives a singular perturbation problem as $\eps \rightarrow 0^{+}$.

Using the method of multiple scales with a continuum slow time variable, this leads to approximate solutions where ${x}(n)$ grows exponentially. The solution breaks down for large $n$ as the dominant balance in the rescaled difference equation changes and a new scaling for $x$ is required. The late time solution based on this new scaling can then be connected with the early time solution using the method of matched asymptotic expansions.
% \begin{subequations}
% \begin{gather}
%  \tilde{x}(n+1) + \tilde{x}(n) 
%  = \frac{2}{9} %\\
%  - \eps \left[ \frac{\tilde{x}(n)}{3} + 3 \,\tilde{x}(n)^2 \right]
%  - \eps^2 \,\tilde{x}(n)^2; 
%  \label{Intro-StaticAlgebraic-Main}
%  \\
%  \tilde{x}(0) = 0,
%  \label{Intro-StaticAlgebraic-IC}
% \end{gather}
% \label{Intro-StaticAlgebraic}
% \end{subequations}

%In Section \ref{S:StaticExpon}, we consider the situation where $\bif = 3+\eps$, but $x(0)$ is chosen to be exponentially close to the unstable 1-periodic equilibrium at $\xoneeqm = \frac{2+\eps}{3+\eps}$. In this case, it takes much longer for the solution trajectory to escape the neighbourhood of $\xoneeqm$. As a result, it becomes necessary to introduce two slow time scales in the method of multiple scales in order to obtain an early time solution that is asymptotically valid up to the point where the dominant balance changes.

%In our analysis of the `exponentially close' problem in section \ref{S:StaticExpon}, we find that it is most convenient to rescale $x(n)$ by considering the size of the perturbation away from the 1-periodic equilibrium, $\xeqm_\text{1per}(n;3+\eps) = \frac{2+\eps}{3+\eps}$, rather than the size of the perturbation away from $\xeqm_\text{1per}(n;3) = \frac{2}{3}$. 

Then, in Section \ref{S:Dynamic}, we demonstrate that our methods can be applied to non-autonomous difference equations by analysing the case where $\bif$ is taken to be a slowly varying function of $n$: $\bif(n) = 3 + \eps^2 \, n$. This introduces a number of subtleties, especially with regard to finding the appropriate asymptotic scalings for $x$ and $n$ in the early time and late time regimes. However, the methods developed in Section \ref{S:StaticAlgebraic} can still be used, and we find that we can obtain a uniformly valid composite asymptotic solution.

%For example, the 1-periodic equilibrium (more properly, the 1-periodic adiabatic manifold) is now a function of $n$, which affects our method for rescaling $x$ in the early time problem. Similarly, it is much more challenging to determine the appropriate rescalings of $x$ and $n$ for the late time problem.

%Additionally, the dependence of $\bif$ on $n$ means that it is less obvious to see the appropriate slow time scale for multiple scales analysis. Also, the rapidity with which $x(n)$ moves away from the unstable 1-periodic manifold introduces complexities to the late time rescalings of $x$ and $n$ in the method of matched asymptotic expansions.

Throughout Sections \ref{S:StaticAlgebraic} and \ref{S:Dynamic}, we focus solely on the discrete logistic equation. Then, in Section \ref{S:Discussion}, we describe how the methods developed in these sections can be applied to a much wider class of difference equations, and we discuss the different combinations of multiple scales and matched asymptotic expansions that could be used as a general strategy for the asymptotic analysis of difference equations. 

%discuss some general features of our results. Section \ref{S:Strategy} discusses how the methods developed in this paper could be used to form the foundation of a general strategy for solving \aside{FINISH THIS}

%As we proceed with our analysis of these three problems, a number of common features emerge, some of which are perhaps unexpected. For example, we find that the rescaling of the time variables in each problem involves translation instead of (or as well as) stretching. This is in contrast to the use of matched asymptotic expansions in `classical' boundary layer problems, where typically it is only necessary to stretch the independent variable. %Additionally, we  find that the composite asymptotic expansion obtained using the method of matched asymptotic expansions is idential to the late time asymptotic expansion, and there are yet other features of our results that suggest a relationship with problems in exponential asymptotics. 
%This and other interesting features of our analysis are discussed further in Section \ref{S:Discussion}. 

\section{The discrete logistic equation with $\bif = 3+\eps$ and $x(0) = \frac{2}{3}$}
\label{S:StaticAlgebraic}

\subsection{Early time scaling}
\label{S:StaticAlgebraicIntro}

In Section \ref{S:StaticAlgebraic}, we will apply the combined method of multiple scales with matched asymptotic expansions to the discrete logistic equation, \eqref{DiscreteLogistic}, with $\bif = 3+\eps$, where $0 < \eps \ll 1$:
\begin{equation}
 x(n+1) = (3+\eps) \, x(n) \, \big[1 - x(n)\big].
 \label{BaseEqn-Static}
\end{equation}
Using \eqref{1PeriodicSol} and \eqref{2PeriodicSol}, 1-periodic and 2-periodic equilibria of \eqref{BaseEqn-Static} are asymptotically given by
\begin{gather}
 \xoneeqm(3+\eps) =  \frac{2+\eps}{3+\eps} = \frac{2}{3} + \frac{\eps}{9} + \mathcal{O}(\eps^2), 
 \label{1PeriodicSol-Eps} \\
 \xtwoeqm^{\pm}(3+\eps) = \frac{2}{3} \pm \frac{\sqrt{\eps}}{3} - \frac{\eps}{18} + \mathcal{O}(\eps^\frac{3}{2}),
 \label{2PeriodicSol-Eps}
\end{gather}
as $\eps \rightarrow 0^{+}$.

% \begin{subequations}
% \begin{gather}
% x(n+1) = (3+\eps) \, x(n) \, \big[1 - x(n)\big],  \\
% x(0) = \frac{2}{3}.
% \end{gather}
% \label{BaseEqn-StaticAlgebraic}
% \end{subequations}

\begin{figure*}
\centering
\includegraphics[]{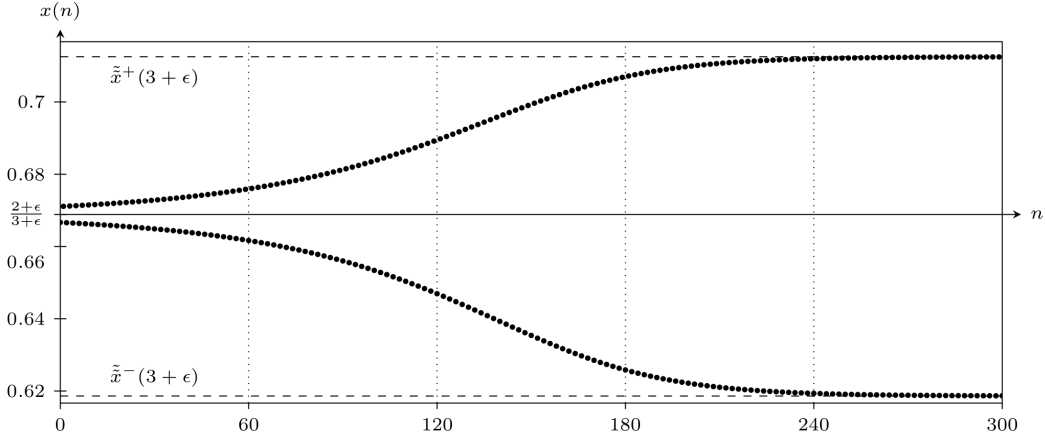}
%\begin{tikzpicture}
%[xscale=0.04*0.9,>=stealth,yscale=1.2*0.8*0.9]
%\draw[->] (0,0) -- (305,0) node[right]{\scriptsize{$n$}};
%\draw[->] (0,2.5) -- (0,2.75) node[above]{\scriptsize{$x(n)$}};
%\draw (1,1.666666667) -- (-2,1.666666667) node[left]{\scriptsize{$0.7$}};
%\draw (2,0.666666667) -- (-2,0.666666667) node[left]{\scriptsize{$0.68$}};
%\draw (2,0) -- (-2,0) node[left]{\scriptsize{$\tfrac{2}{3}$}};
%\draw (2,-0.33333333) -- (-2,-0.33333333) node[below left]{\scriptsize{$0.66$}};
%\draw (2,-1.33333333) -- (-2,-1.33333333) node[left]{\scriptsize{$0.64$}};
%\draw (2,-2.33333333) -- (-2,-2.33333333) node[left]{\scriptsize{$0.62$}};
%\draw [dotted] (60,-2.5) -- (60,2.5);
%\draw [dotted] (120,-2.5) -- (120,2.5);
%\draw [dotted] (180,-2.5) -- (180,2.5);
%\draw [dotted] (240,-2.5) -- (240,2.5);
%\draw (0,-2.5) -- (0,-2.6) node[below] {\scriptsize{$0$}};
%\draw (60,-2.5) -- (60,-2.6) node[below] {\scriptsize{$60$}};
%\draw (180,-2.5) -- (180,-2.6) node[below] {\scriptsize{$180$}};
%\draw (120,-2.5) -- (120,-2.6) node[below] {\scriptsize{$120$}};
%\draw (240,-2.5) -- (240,-2.6) node[below] {\scriptsize{$240$}};
%\draw (300,-2.5) -- (300,-2.6) node[below] {\scriptsize{$300$}};
%\begin{scope}
%\clip (0,-2.5) -- (0,2.5) -- (300,2.5) -- (300,-2.5) -- cycle;
%\draw[black] plot[only marks, mark=*,mark options={fill=black,xscale= 25*0.45*0.8,yscale=0.8333*0.45}] file {xvec_Fig2.txt};
%\end{scope}
%\draw (0,-2.5) -- (0,2.5) -- (300,2.5) -- (300,-2.5) -- cycle;
%\end{tikzpicture}
\caption{Numerical solution to the discrete logistic equation in the neighbourhood of its first period doubling bifurcation. Black circles indicate the solution $x(n)$ to \eqref{BaseEqn-Static} where $\eps = 0.02$ and $x(0) = \frac{2}{3}$. As $n$ increases, we see that $x(n)$ moves away from the unstable equilibrium at $\xoneeqm = \frac{2+\eps}{3+\eps} \approx 0.6689$, before settling into the oscillatory pattern given by equation \eqref{2PeriodicSol-Eps}, with the period-two manifold illustrated by dashed curves.}
\label{Static_exact}
\end{figure*}

A numerical solution of \eqref{BaseEqn-Static} where the initial condition is close to $\xoneeqm$ is presented in Figure \ref{Static_exact}. We note that the solution gives the appearance of a smooth envelope around a system of oscillations away from the unstable steady state at $\xoneeqm$ towards the 2-periodic equilibrium given by \eqref{2PeriodicSol-Eps}.

Consider the case where we wish to solve \eqref{BaseEqn-Static} subject to the initial condition $x(0) = \frac{2}{3}$. As given in \eqref{BaseEqn-Static}, this difference equation is strongly nonlinear. However, by introducing the rescaling
\begin{equation}
 x(n) = \tfrac{2}{3} + \eps\, X(n)
 \label{StaticEarlyRescaling}
\end{equation}
we obtain the weakly nonlinear problem,
%\begin{subequations}
\begin{equation}
 X(n+1) + X(n)  = \frac{2}{9} 
 - \eps \left[ \frac{X(n)}{3} + 3 \,X(n)^2 \right]
 - \eps^2 \,X(n)^2,
 \label{StaticAlgebraic-EarlyMain}
\end{equation}
with $X(0)=0$.
%\label{StaticAlgebraic-Early}
%\end{subequations}

%This problem is well-suited to being solved using the method of multiple scales. However, we note from the start that there is a strong possibility that we will run into difficulties. Since the 2-periodic equilibrium stated in \eqref{2PeriodicSol-Eps} is $\mathcal{O}(\eps^{\frac{1}{2}})$ away from $x = \frac{2}{3}$, it is likely that the scaling we have introduced in \eqref{StaticEarlyRescaling} will fail as we leave the neighbourhood of the 1-periodic equilibrium and approach the 2-periodic equilibrium.

%Throughout the rest of this Section, we will move relatively slowly through the manipulations associated with the method of multiple scales and then the method of matched asymptotic expansions. Then, having established the method with this example, we will move more quickly in sections \ref{S:StaticExpon} and \ref{S:Dynamic}.

\subsection{Early time solution}
\label{S:StaticAlgebraicEarlySol}

We begin by introducing $t = \eps \, n$ as a slow time variable alongside the fast time variable, $n$, and using the method of multiple scales where $t$ is treated as a continuum variable.
%
% First, consider what would happen if we were to propose a simple asymptotic expansion for $\tilde{x}(n)$ of the form
% \begin{equation}
%  \tilde{x}(n) \sim \tilde{x}_0(n) + \eps \, \tilde{x}_1(n) + \ldots
%  \label{StaticNaive}
% \end{equation}
% and substitute into \eqref{StaticAlgebraic-Early}. Solving the difference equations obtained by collecting terms of $\mathcal{O}(1)$ and $\mathcal{O}(\eps)$ in \eqref{StaticAlgebraic-Early}, we would find that
% \begin{gather}
%  \tilde{x}_0(n) = \frac{1 + (-1)^n}{9}, \\
%  \tilde{x}_1(n) = \frac{-1+ (-1)^n}{18} - \frac{(-1)^n}{9} \, n.
% \end{gather}
% 
% This is a classic example of a secular term appearing in an asymptotic expansion. Our solution based on \eqref{StaticNaive} fails to be asymptotic when $n = \mathcal{O}(\eps^{-1})$ because the correction term $\tilde{x}_1(n)$ grows to be the same size as the leading order solution, $\tilde{x}_0(n)$. Thus, it is appropriate for us to introduce $t = \eps n$ as a slow time variable alongside the fast time variable $n$, and use the method of multiple scales to remove the secularity. In this section (as throughout this paper), $t$ will be treated as a continuum variable.
%
Our multiple scales ansatz takes the form
$
 X(n) \equiv X(n,t)$,
where $X(n,t)$ will be expanded as an asymptotic series in powers of $\eps$. 

Following Hoppensteadt and Miranker \cite{Hoppensteadt1977}, $X(n+1)$ can be expanded as a Taylor series of the form
\begin{equation}
X(n+1,t+\eps) = \sum_{j=0}^{\inf} \frac{\eps^j}{j!} \diff{^j X(n+1,t)}{t^j}.
\end{equation} 
Substituting into \eqref{StaticAlgebraic-EarlyMain}, we therefore find that the equation to solve is
\begin{align}
%\sum_{j=0}^{\inf} \frac{\eps^j}{j!} \diff{^j X(n+1,t)}{t^j}
 \nonumber X(n+1,t) + X(n,t) =  \frac{2}{9} 
- \eps&\left(\frac{X(n,t)}{3} + 3X^2(n,t) +\diff{X(n+1,t)}{t} \right) \\
&- \eps^2\left(X^2(n,t) + \frac{1}{2} \diff{^2X(n+1,t)}{t^2}\right) + \mathcal{O}(\eps^3),\label{1.0 ScaledStatic2}
\end{align}
subject to the condition that $X(0,0) = 0$. This is now a \DeltaDE{} expressed in terms of the discrete variable $n$ and the continuous variable $t$. 

We expand $X(n,t)$ as an asymptotic series in the limit $\eps \rightarrow 0^+$:
\begin{equation}\label{1.0 Series}
X(n,t) \sim \sum_{k=0}^{\inf}\eps^k X_k(n,t).
\end{equation}
Applying this series to (\ref{1.0 ScaledStatic2}), we collect terms of the same order in $\eps$ to obtain a system of difference equations for the functions $X(n,t)$. Collecting terms of $\Ord(1)$ gives% and $\Ord(\eps)$ gives
\begin{equation}\label{1.0 FirstOrder}
 \mathcal{O}(1):  \quad 
 X_{0}(n+1,t) + X_{0}(n,t) = \tfrac{2}{9},%\qquad &X^{(0)}(0,0) &= 0,\\
\end{equation}
%\begin{multline}
 %\label{1.0 SecondOrder1}
 %\mathcal{O}(\eps): \quad  
 %X_{1}(n+1,t) + X_{1}(n,t) = -\frac{1}{3}X_{0}(n,t) \\
 %- 3\left[X_{0}(n,t)\right]^2 - \diff{X_{0}(n+1,t)}{t},%\qquad &X^{(1)}(0,0) &= 0.
%\end{multline}
with $X_{0}(0,0) = 0$;

Solving (\ref{1.0 FirstOrder}), we find that
\begin{equation}\label{1.0 Xn1}
X_{0}(n,t) = \tfrac{1}{9}[1-(-1)^n A(t)],\qquad A(0) = 1.
\end{equation}
where $A(t)$ is some arbitrary smooth function of $t$. Collecting terms at $\mathcal{O}(\eps)$ gives an equation for the next term in the series, 
\begin{equation}\label{1.0 SecondOrder}
X_{1}(n+1,t) + X_{1}(n,t) 
=  - \tfrac{2}{27} - \tfrac{1}{27}\left[A(t)\right]^2 - \tfrac{1}{9}(-1)^n(A'(t) - A(t)).
\end{equation}

%Noting that $(-1)^n$ is a solution to the homogeneous equation associated with \eqref{1.0 SecondOrder}, we find that $A'(t)-A(t) = 0$ in order for our solution to \eqref{1.0 SecondOrder} to be free of secular terms in $n$. 
In order to remove the secular $(-1)^n$ terms from the right hand side of \eqref{1.0 SecondOrder}, as is standard in the method of multiple scales, we impose $A'(t) - A(t) = 0$.
Combined with the boundary condition in \eqref{1.0 Xn1}, this gives $A(t) = \e^{t}$, and hence $X_{0}(n,t) = \tfrac{1}{9}[1 - (-1)^{n} \, \e^t]$.
%\begin{equation}\label{1.0 x_n 1term}
 %X_{0}(n,t) = \tfrac{1}{9}[1 - (-1)^{n} \, \e^t].
%\end{equation}

% Having eliminated the $(-1)^n$ terms from its right hand side, equation \eqref{1.0 SecondOrder} can now be solved to obtain $X_1(n,t)$ as
% \begin{equation}
% X_{1}(n,t) = \frac{[1-\e^{2t}][1-(-1)^n]}{54} + \frac{(-1)^n B(t)-1}{18},%\qquad B(0) = 1,
% \label{EarlyX1-Temp}
% \end{equation}
% where $B(t)$ is again an arbitrary smooth function of $t$, and $B(0) = 1$ in order to satisfy the initial condition $X_1(0,0) = 0$. Going to the $\Ord(\eps^2)$ terms in \eqref{1.0 GeneralTerms}, we find that $B(t)$ must satisfy
% \begin{equation}
% B'(t) - B(t) = -\frac{1}{3}+\e^t-\frac{\e^{2t}}{3} + \frac{2\e^{3t}}{9}% \qquad B(0) = 1,
% \label{EarlyBtEquation}
% \end{equation}
% in order to avoid secularity in $n$, and hence \eqref{EarlyX1-Temp} becomes
% Repeating the same process at the next order, we find that
% \begin{equation}\label{1.0 X1}
% X_1(n,t) = \frac{6+3\e^{2t}-(-1)^n\left(8+9t+\e^{2t}\right)\e^t}{162}.
% \end{equation}
% Thus, the two term solution for $X(n,t)$ is
% \begin{multline}\label{1.0 Xn_2terms}
% X(n,t) = \frac{1-(-1)^n\e^t}{9} \\
% - \eps\left[\frac{6+3\e^{2t}-(-1)^n\left(8+9t+\e^{2t}\right)\e^t}{162} \right] + \mathcal{O}(\eps^2).
% \end{multline}

As described in Appendix \ref{static_appendix}, this process can be repeated methodically, leading to the result that the solutions $X_r(n,t)$  may be written in the form
\begin{equation}\label{1.0 Xr}
X_{r}(n,t) = f_r(t) + g_r(t) \, (-1)^{n},
\end{equation}
where $f_r(t)$ is given in \eqref{1.0 fk soln} and $g_r(t)$ is obtained by solving \eqref{1.0 gk ode}.

%are given by
%\begin{multline}
% f_r(t) = 
%  - \frac{1}{6} f_{r-1}(t) \\
%  - \frac{3}{2} \sum_{k=0}^{r-1} \Big[ f_{k}(t) \, f_{r-1-k}(t) + g_{k}(t) \, g_{r-1-k}(t) \Big] \\
%  - \frac{1}{2} \sum_{k=0}^{r-2} \Big[ f_k(t) \, f_{r-2-k}(t) + g_k(t) \, g_{r-2-k}(t) \Big] \\
%  - \frac{1}{2} \sum_{j=1}^{r}\frac{1}{j!} \, f_{r-j}^{(j)}(t), \label{StaticEarlyStanding}
%\end{multline}
%and the functions $g_r(t)$ are obtained by solving
%\begin{multline}
% g_r'(t) -  g_r(t)
%  =  6 \sum_{k=1}^{r}  f_{k}(t) \, g_{r-k}(t) \\   
%  + 2 \sum_{k=0}^{r-1} f_k(t) \, g_{r-1-k}(t) 
%  - \sum_{j=1}^{r}\frac{1}{j!} \,  g_{r-j}^{(j+1)}(t), \label{StaticEarlyOscillatingDE}
%\end{multline}
%subject to the initial condition $g_r(0) = -f_r(0)$. %Again, sums from 0 to $-1$ are taken to be zero.
%
%Applying these formulae to the case where $r=1$,  we find that 
%\begin{equation}\label{1.0 X1}
%X_1(n,t) = \tfrac{1}{162}\left[6+3\,\e^{2t}-(-1)^n\,\left(8\,\e^t+9\,t\,\e^t+\e^{3t}\right)\right],
%\end{equation}
%and hence
Using these expressions to calculate $X_1(n,t)$, we find that
\begin{equation}\label{1.0 Xn_2terms}
X(n,t) = \tfrac{1}{9}[{1-(-1)^n \, \e^t}] 
- \tfrac{\eps}{162}\left[{6+3\,\e^{2t}-(-1)^n\,\left(8\,\e^t+9\,t\,\e^t+\e^{3t}\right)} \right] + \mathcal{O}(\eps^2).
\end{equation}

\begin{figure*}
\centering
\subfloat[One-term approximation]{
\includegraphics[]{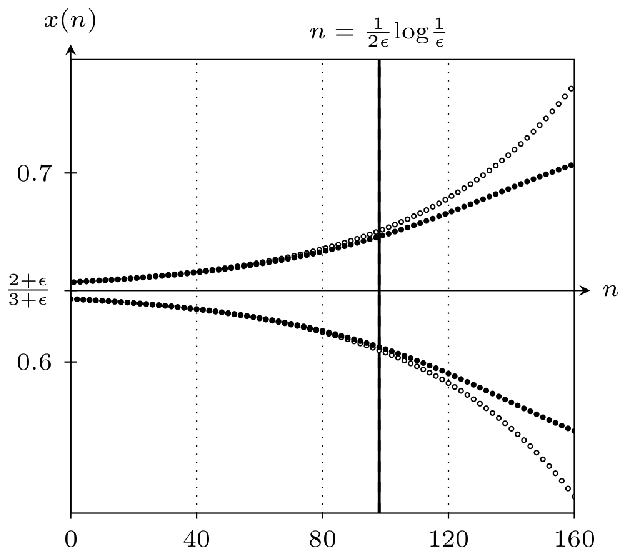}
\label{F:OneTermStaticEarly}
}
\subfloat[Two-term approximation]{
\includegraphics[]{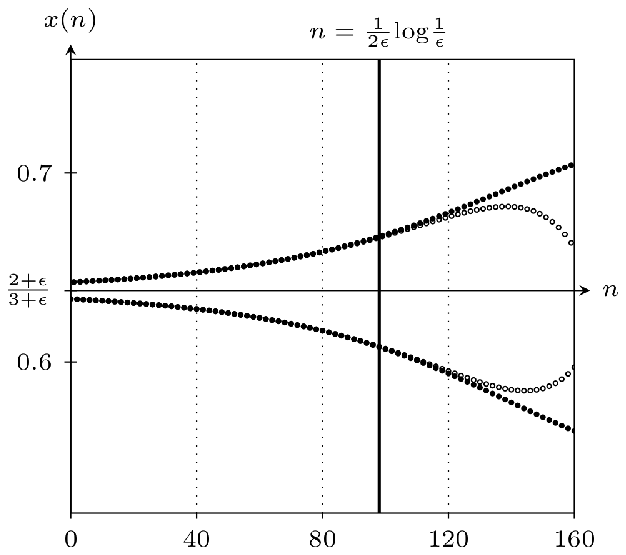}
\label{F:TwoTermStaticEarly}
}
\caption{Numerical solution $x(n)$ to \eqref{StaticAlgebraic-EarlyMain} with $x(0) = \frac{2}{3}$ compared with (a) the one-term early-time approximation  and (b) the two-term early-time approximation, both obtained from (\ref{1.0 Xn_2terms}), in the case where $\eps = 0.02$. In each figure, the exact solution is represented by black circles, while the asymptotic approximation is represented by white circles. The point $n = \tfrac{1}{2\eps}\mathrm{log}\tfrac{1}{\eps}$ is depicted as a thick line; this represents the approximate time at which the series ceases to be asymptotic, as discussed in Section \ref{S:StaticAlgebraicFailure}. In both cases, the approximation fails when continued beyond this point. \label{F:StaticEarly}} 
\end{figure*}

Figure \ref{F:StaticEarly} compares a numerical solution of \eqref{StaticAlgebraic-EarlyMain} with asymptotic solutions based on the one term approximation %\eqref{1.0 x_n 1term} %(shown in Figure \ref{F:OneTermStaticEarly}) 
and the two term approximation, obtained by truncating \eqref{1.0 Xn_2terms}. %(shown in Figure \ref{F:TwoTermStaticEarly}). 
We see that in both cases the approximation obtained from \eqref{1.0 Xn_2terms} accurately describes the solution trajectory in its early stages; eventually, however, both asymptotic series approximations become inaccurate. Furthermore, using \eqref{1.0 Xr} to generate subsequent terms in the series does not improve the accuracy beyond this point.

\subsection{Failure of the early time expansion}
\label{S:StaticAlgebraicFailure}

The method of multiple scales is useful for considering the slow accumulation of small deviations; however, for changes in the dominant balance of an equation, the method of matched asymptotics is necessary \cite{Lagerstrom-MAE}.

This is not straightforward in the present problem. Considering the asymptotic sizes of the various terms in \eqref{StaticAlgebraic-EarlyMain} as $X \rightarrow \infty$, we note that the $X(n+1)$ and $X(n)$ terms will always dominate the $\frac{2}{9}$ term, and that the $-3 \, \eps \, X(n)^2$ term will dominate all of the other terms. Hence, a na{\"{\i}}ve inspection of \eqref{StaticAlgebraic-EarlyMain} would suggest that a new dominant balance is attained when $X = \ord\big[\eps^{-1}\big]$, which is incorrect. In fact, subsequent analysis in this section will show that the failure of the asymptotic series occurs earlier, when $X = \ord\big[\eps^{-\frac{1}{2}}\big]$.

This error occurs because the failure of the early time asymptotic expansion is not associated with the discrete-scale behaviour of solutions to \eqref{StaticAlgebraic-EarlyMain}, but rather due to the exponential growth of the continuous slow time envelope. This growth causes the differential equations associated with the multiple scales formulation to undergo a change in dominant balance. To determine when this occurs, we can ignore discrete-scale variation, and reformulate \eqref{StaticAlgebraic-EarlyMain} purely in terms of the slow time variable, $t$.

%In fact, the failure occurs because the continuous envelope grows exponentially and causes the differential equations associated with the multiple scales formulation to undergo a change in dominant balance. 

The most natural way is to remove the effects of the discrete scale by considering the doubled map, as it is clear from the early-time analysis that the behaviour on the discrete scale is 2-periodic. Applying \eqref{StaticAlgebraic-EarlyMain} twice, we find that
\begin{align}
 \nonumber X(n+2) - X(n) 
 = &\quad \eps \, \left[2 \, X(n) - \tfrac{2}{9}\right] 
 + \eps^2 \, \left[-18 \, X(n)^3+ 3 \, X(n)^2 + X(n) - \tfrac{4}{81} \right] \\
 \nonumber &+ \eps^3 \, \left[-27 \, X(n)^4 - 18 \, X(n)^3 + 2 \, X(n)^2 + \tfrac{4}{27}X(n) \right] \\
 \nonumber &+ \eps^4 \, \left[-27 \, X(n)^4 - 6 \, X(n)^3 + \tfrac{1}{3}X(n)^2\right] \\
 &- \eps^5 \, \left[9 \, X(n)^4 + \tfrac{2}{3}X(n)^3\right]
 - \eps^6 \, X(n)^4.
 \label{1.0 DoubledMap}
\end{align} 

%If we then introduce the doubled discrete time variable, $\tilde{n}$, so that $n = 2 \, \tilde{n}$, and apply the method of multiple scales to \eqref{1.0 DoubledMap} with fast time variable $\tilde{n}$ and continuum slow time variable $t = 2 \, \eps \, \tilde{n}$, we would quickly find that $X$ is independent of $\tilde{n}$. %This is because the homogeneous equation equivalent to the right hand side of \eqref{1.0 DoubledMap} is solved by $(1)^{\tilde{n}} \equiv 1$, and this is the only form in which $\tilde{n}$ can appear.

Thus, we could skip the method of multiple scales entirely by proposing the continuum ansatz $X(n) \equiv \chi(t)$ in \eqref{1.0 DoubledMap} with $t = \eps \, n$, leading to the infinite-order ODE
\begin{align}
 \nonumber 2  \, \underset{\textrm{\ding{172}}}{\diff{\chi}{t}} + \sum_{j=2}^{\inf} \frac{\eps^{j-1} \, 2^j}{j!} \diff{^j \chi}{t^j} 
 =&   \left[\underset{\textrm{\ding{173}}}{2 \, \chi(t)} - \tfrac{2}{9}\right] 
 + \eps \, \left[\underset{\textrm{\ding{174}}}{-18 \, \chi(t)^3} + 3 \, \chi(t)^2 + \chi(t) - \tfrac{4}{81} \right] \\
 \nonumber & + \eps^2 \, \left[-27 \, \chi(t)^4 - 18 \, \chi(t)^3 + 2 \, \chi(t)^2 + \tfrac{4}{27} \chi(t) \right] \\
 \nonumber &+ \eps^3 \, \left[-27 \, \chi(t)^4 - 6 \, \chi(t)^3 + \tfrac{1}{3} \chi(t)^2\right] \\
 &- \eps^4 \, \left[9 \, \chi(t)^4 + \tfrac{2}{3} \chi(t)^3\right]
 - \eps^5 \, \chi(t)^4.
 \label{1.0 DoubledMapContinuum}
\end{align}

This is now in a form where we can apply the usual methods for finding the new dominant balance; we can propose rescalings for $\chi$ and $t$ and seek a new dominant balance in \eqref{1.0 DoubledMapContinuum} that is consistent with long time behaviour of the leading order solution \eqref{1.0 Xn1}. Note that the exponential growth of the leading order solution $X_0(n,t)$ means that it will be necessary to propose an affine rescaling of $t$ rather than a linear rescaling.

%which upon rescaling gives $\chi_0(t) = \ord(\e^t)$, where $\ord$ is used to represent `strict order' as in \cite{HinchPert}, so that $f = \ord(g)$ is equivalent to the case where $f = \Ord(g)$ but $f \neq o(g)$.  

Specifically, we make the rescalings $\chi = \delta \, \xi$ and $ t= K_0 + K_1 \, s$, where $\xi$ and $s$ are $\Ord(1)$ variables, and where $\delta$, $K_0$, and $K_1$ depend only on $\eps$ with $\delta \gg 1$, $K_0 \gg 1$ and $K_0 \gg K_1$ as $\eps \rightarrow 0^+$. Substituting into the leading order solution from \eqref{1.0 Xn_2terms}, we find that
\begin{equation}
 \delta \, \xi = \ord\left[ \e^{K_0} \, \e^{K_1 \, s}\right],
 \label{PreparationForRebalancing-Static}
\end{equation}
where $\ord$ is used to represent `strict order' as in \cite{HinchPert}, so that $f = \ord(g)$ is equivalent to the case where $f = \Ord(g)$ but $f \neq o(g)$. 

Since \eqref{PreparationForRebalancing-Static} must hold whenever $\xi$ and $s$ are both $\ord(1)$, we find that it is appropriate to choose $K_1 = 1$ and $K_0 = \log \delta$. Substituting into the doubled map \eqref{1.0 DoubledMapContinuum}, we find that the terms labelled \ding{172} and \ding{173} always remain part of the dominant balance, and the first term to grow to be the same size as \ding{172} and \ding{173} is the one labelled \ding{174}. This occurs when $\eps \, \delta^3 = \ord(\delta)$; hence, it is appropriate to take $\delta = \eps^{-\frac{1}{2}}$ and $K_0 = \frac{1}{2}\log\left(\frac{1}{\eps}\right)$. This suggests that the asymptotic series breaks down due to a change in dominant balance when $X = \ord[\eps^{-\frac{1}{2}}]$, which occurs at $t = \eps\, n = \frac{1}{2}\log\left(\frac{1}{\eps}\right) + \Ord(1)$.

\subsection{Late time solution and matched aysmptotic expansions}
\label{S:StaticAlgebraicLate}

Using the doubled map, we determined the scaling of $X$ and $t$, and hence $x$ and $n$, at which the early-time expansion loses asymptoticity. Hence, we apply the rescalings
\begin{align}
 x(n) &= \tfrac{2}{3} + \eps^{1/2}\xi(m,s), \label{StatLate1} \\
 t &= \tfrac{1}{2}\log(\tfrac{1}{\eps}) + s, \label{StatLate2} \\
 n &= \tfrac{1}{2\eps}\log(\tfrac{1}{\eps}) -\gamma + m, \label{StatLate3}
\end{align}
 %which led to the rescalings given in equations \eqref{1.0 Xxi Rescaling}, \eqref{1.0 ts Rescaling}, and \eqref{1.0 mn Rescaling}. 
where $\gamma$ is a constant chosen so that $0 \leq \gamma < 2$, and $m-n \equiv 1 \pmod 2$. The latter condition is chosen for convenience in subsequent analysis. Applying these rescalings to \eqref{StaticAlgebraic-EarlyMain} and using the method of multiple scales with discrete fast time variable $m$ and continuum slow time variable $s$, we recover
\begin{multline} \label{1.1 OuterGoverning}
\xi(m+1,s) + \xi(m,s) 
= \eps^{\frac{1}{2}}\, \left[\frac{2}{9} - 3\,\xi(m,s)^2\right] \\
- \eps\left[\tfrac{1}{3} \,\xi(m,s) + \diff{\xi(m+1,s)}{s}\right] - \eps^{\frac{3}{2}}\, \xi(m,s)^2 + \mathcal{O}(\eps^2)
\end{multline}

Since \eqref{1.1 OuterGoverning} involves powers of $\eps^\frac{1}{2}$, we introduce a series expansion for $\xi(m,s)$ in half-powers of $\eps$ of the form
\begin{equation}
\xi(m,s) \sim \sum_{k=0}^{\inf} \eps^{\frac{k}{2}} \, \xi_{k}(m,s).
\end{equation}

Collecting terms of the same order in $\eps$, we obtain
\begin{align}
  \label{1.1 O1}
  \mathcal{O}(1):& \quad  
  \xi_0(m+1,s) + \xi_0(m,s) = 0,\\
  \mathcal{O}(\eps^{\frac{1}{2}}):&  \quad  
  \xi_1(m+1,s) + \xi_1(m,s) = \tfrac{2}{9} - 3 \, \xi_0(m,s)^2, 
\end{align}
%\begin{multline}
  %\mathcal{O}(\eps): \quad  
  %\xi_2(m+1,s) + \xi_2(m,s) = -\frac{1}{3} \, \xi_0(m,s) \\
  %-6\, \xi_0(m,s)\,\xi_1(m,s)
  %- \diff{\xi_0(m+1,s)}{s},
%\end{multline}
and so on.

%in general for $r \geq 2$,
%\begin{multline}
%\mathcal{O}(\eps^{\frac{r}{2}}):  \quad 
%\xi_r(m+1,s)  + \xi_r(m,s) 
%= -\frac{1}{3} \, \xi_{r-2}(m,s) \\
%-3\sum_{k=1}^{r-1} \xi_k(m,s) \, \xi_{r-k}(m,s)
%-\sum_{k=0}^{r-3} \xi_{k}(m,s) \, \xi_{r-k}(m,s) \\
%- \sum_{j=1}^{\lfloor\frac{r}{2}\rfloor}\frac{1}{j!} \, \diff{^j \xi_{r-2j}(m+1,s)}{s^j}.
%\end{multline}

At leading order, we can solve (\ref{1.1 O1}) easily to find that $\xi^{(0)}(m,s) = (-1)^m P(s)$, where $P(s)$ is an arbitrary smoothly-varying function in $s$. At $\mathcal{O}(\eps^{\frac{1}{2}})$, we do not obtain a secularity condition. At the next order, however, we find the secularity condition for $P(s)$ is given by $P'(s) = P(s)-9\,P^3(s)$. Solving this yields
\begin{equation}
P(s) = \frac{\e^{s}}{3\sqrt{\e^{2s}+ \kappa_p}}, \label{1.0 P}
\end{equation}
where $\kappa_p$ is an arbitrary constant, and we recall that $\xi_0 = (-1)^m P(s)$.

%but find that
%\begin{equation}
%\xi_1(m,s)  = \tfrac{1}{18}{[2-27\, P(s)^2] \, [1-(-1)^m]} + (-1)^m \, Q(s), \label{1.0 xi1}
%\end{equation}
%where $Q(s)$ is again an arbitrary smoothly-varying function. At the next order, however, we find that
%\begin{multline}
%\xi_2(m+1,s) + \xi_2(m,s) = -6\, P(s)\,Q(s)+\frac{2}{3} \, P(s)-9\,P(s)^3 \\
%+ (-1)^m\left[-P(s)+9 \, P(s)^3+P'(s)\right].
%\end{multline}
%This equation gives%hence
%\begin{equation}
% \xi_0(m,s) = \frac{(-1)^m \, \e^{s}}{3\sqrt{\e^{2s}+ \kappa_p}}, \label{1.0 xi0}
%\end{equation} 

In order to determine the value of the constant $\kappa_p$, we apply Van Dyke's matching criterion (see, for example, \cite{HinchPert}) to match the behaviour of the solution in the late-time region with the known behaviour in the early-time region. %Equivalently, we could introduce an intermediate time variable, $T$, so that $t = \alpha(\eps) + T$ where $1 \ll \alpha \ll \frac{1}{2} \log(\frac{1}{\eps})$ and use intermediate variable matching; Van Dyke matching is sufficient, however, and using Van Dyke makes it straight forward to obtain the uniformly valid composite solution. 

%At leading order, Van Dyke's matching criterion involves expanding $X^{(0)}(n,t)$ in terms of late-time variables, expanding $\xi_0(m,s)$ in terms of early-time variables, and equating the two at leading order.
Van Dyke's matching criterion involves expanding the early time solution in late time variables and equating this with an appropriate expansion of the early time solution in late time variables. To determine $\kappa_p$, we express the one-term late time solution $\xi_0(m,s)$ in terms of early time variables $n$ and $t$, and equate the leading-order behaviour of this expression with the leading-order behaviour of the one-term early time solution $X_0(n,t)$ in terms of late time variables $m$ and $s$. This gives the matching condition
%
%we find that
%\begin{align}
% \xi_0 &= (-1)^{n+1} \, \frac{\exp [t - \tfrac{1}{2} \log(\tfrac{1}{\eps})]}{3\sqrt{\exp[2\,t - \log(\tfrac{1}{\eps})] + \kappa_p}} \notag \\
% &= (-1)^{n+1} \,\eps^{\frac{1}{2}} \, \frac{\e^{t}}{3 \, \sqrt{\eps \, \e^{2t} + \kappa_p}}.
%\end{align}
%Taking one term of this expression as $\eps \rightarrow 0$, we find that 
%\begin{equation}
% \eps^{-\frac{1}{2}} \xi_0 \sim (-1)^{n+1} \, \frac{\e^t}{3 \, \sqrt{\kappa_p}}, \label{FirstMatchA}
%\end{equation}
%which gives us the `matching term' based on matching the one-term late time solution with the one-term early time solution. Van Dyke's matching principle requires that this be equivalent to the leading-order behaviour of the one-term early time solution expressed in late time variables.
%Taking the one-term early time solution \eqref{1.0 Xn1} and expressing it in terms of late time variables, we find that
%\begin{equation}
% X_0 = (-1)^{m} \, \frac{\e^s}{9} + \Ord[\eps \, \log (\tfrac{1}{\eps}) ] \label{FirstMatchB}
%\end{equation}
%To match between \eqref{FirstMatchA} and \eqref{FirstMatchB}, we therefore require that
\begin{equation}
\frac{(-1)^m\,\e^{s}}{9\, \sqrt{\eps}} \equiv -\frac{(-1)^n \, \e^{t}}{3 \,\sqrt{\kappa_p}},
\end{equation}
which, using \eqref{StatLate1}-\eqref{StatLate3} and exploiting the simplifying assumption that $m -n \equiv 1 \pmod 2$, gives $\kappa_p = 9$.

%The one-term late time solution given by \eqref{1.0 xi0} with $\kappa_p = 9$ is only accurate to $\Ord(\eps^{\frac{1}{2}})$, whereas the one-term early time solution, \eqref{1.0 x_n 1term}, is accurate to $\Ord(\eps)$. In order to construct a late time solution that is also accurate to $\Ord(\eps)$, and hence build a composite solution that is uniformly valid to $\Ord(\eps)$, we need to obtain $\xi_1(m,s)$.  

By continuing to match powers of $\eps$, we obtain a form for $\xi_1(m,s)$, with an appropriate secularity condition. Van Dyke's matching condition can also be used at higher orders, and hence we are able to determine subsequent terms in the asymptotic series, with the two-term expansion given by
\begin{equation}\label{1.0 late-time}
\xi(m,s) = (-1)^m \, \frac{\e^{s}}{3\,\sqrt{9+\e^{2s}}}
+ \frac{18 - \e^{2s}}{162 + 18 \, \e^{2s}} \, \eps^{\frac{1}{2}} 
+ \mathcal{O}(\eps).
\end{equation}

%Continuing to the next order in \eqref{1.1 OuterGoverning} leads to the secularity condition for $Q(s)$, 
%\begin{equation}
% Q'(s) = -\frac{162+9\,\e^{2s}+2\,\e^{4s} }{18\,(9+\e^{2s})^2}- \frac{(9-2\,\e^{2s})\,Q(s)}{9+\e^{2s}},
%\end{equation}
%which gives
%\begin{equation}
% Q(s) = \frac{18-\e^{2s}}{18\,(9+\e^{2s})} + \frac{\kappa_q \,\e^{s}}{(9+\e^{2s})^{\frac{3}{2}}}. \label{1.0 Q}
%\end{equation}

%Comparing the two-term early time solution from \eqref{1.0 Xn_2terms} with the two-term late time solution ultimately yields the result that $\kappa_q = 0$. Substituting back, we find that the two-term expansion at late time is given by

Since the one-term early time approximation from \eqref{1.0 Xn_2terms}, and the two-term late time approximation, \eqref{1.0 late-time}, are both accurate up to $\Ord(\eps)$ in $x$, we can combine them to obtain a composite asymptotic approximation that will also be uniformly valid to $\Ord(\eps)$. This is obtained by converting both approximations into equivalent variables, adding them together, and subtracting the matching term associated with Van Dyke's criterion. In this case, we find that the matching term is identical to the one-term early time approximation, and hence we recover the result that \eqref{1.0 late-time} is also the uniformly-valid composite approximation.

In classical boundary layer problems, it is unusual for the uniformly-valid composite solution to be completely identical to one of the two solutions involved in the matching, as is the case here. It is important to note that this does not imply that it would have been possible to obtain the uniformly-valid solution without considering the early time and late time cases separately. 

Knowing that $X(0) = 0$, we see that the initial condition in this case should take the form $\xi(0,0) = 0$. However, if we take the late time solution from \eqref{1.0 late-time} and consider what happens when $n = 0$, and hence $s = -\frac{1}{2}\log(\frac{1}{\eps})$ and $(-1)^m =-1$, we find that \eqref{1.0 late-time} becomes
\begin{equation}
 \xi_\text{initial} = -\frac{\sqrt{\eps}}{3\,\sqrt{9+\eps}} \cdot 1 + \frac{18 - \eps}{162 + 18\, \eps} \cdot \eps^{\frac{1}{2}} + \text{``}\Ord(\eps)\text{''}.
\end{equation}
This is neither a valid asymptotic expansion (since both terms are the same size), nor an obvious way of representing the true initial condition of $\xi(0,0) = 0$. While the two-term late time asymptotic expansion is \emph{equivalent} to the composite approximation that is uniformly valid up to $\Ord(\eps)$, it is not a Poincar{\'e} expansion when  $s = -\frac{1}{2}\log(\frac{1}{\eps}) + \Ord(1)$. As a result, there is no power series ansatz by which \eqref{1.0 late-time} could be obtained directly from the initial condition; we need to use the method of matched asymptotic expansions to obtain the late time solution based on matching with the early time solution, and then combine the two to obtain a composite approximation.

%Demonstrations of these results (and a generalisation to the other cases considered in this paper) are given in Appendix \ref{A:LateTimeComposite}.

%In this case, however, it is possible to show that matching the $r$-term early time expansion with the $R$-term late time expansion will always yield a situation where the late time solution is identical to the composite approximation as long as $R \geq 2\,r$. 

\begin{figure*}
\centering
%\begin{tikzpicture}
%[xscale=0.04*0.9,>=stealth,yscale=1.2*0.8*0.9]
%\draw (0,-2.5) -- (0,2.5) -- (250,2.5) -- (250,-2.5) -- cycle;
%\draw[->] (0,0) -- (255,0) node[right]{\scriptsize{$n$}};
%\draw[->] (0,2.5) -- (0,2.75) node[above]{\scriptsize{$x$}};
%\draw (2,2) -- (-2,2) node[left]{\scriptsize{$2$}};
%\draw (2,1) -- (-2,1) node[left]{\scriptsize{$1$}};
%\draw (2,0) -- (-2,0) node[left]{\scriptsize{$0$}};
%\draw (2,-1) -- (-2,-1) node[left]{\scriptsize{$-1$}};
%\draw (2,-2) -- (-2,-2) node[left]{\scriptsize{$-2$}};
%\draw [dotted] (60,-2.5) -- (60,2.5);
%\draw [dotted] (120,-2.5) -- (120,2.5);
%\draw [dotted] (180,-2.5) -- (180,2.5);
%\draw [dotted] (240,-2.5) -- (240,2.5);
%\draw (0,-2.5) -- (0,-2.6) node[below] {\scriptsize{$0$}};
%\draw (60,-2.5) -- (60,-2.6) node[below] {\scriptsize{$60$}};
%\draw (180,-2.5) -- (180,-2.6) node[below] {\scriptsize{$180$}};
%\draw (120,-2.5) -- (120,-2.6) node[below] {\scriptsize{$120$}};
%\draw (240,-2.5) -- (240,-2.6) node[below] {\scriptsize{$240$}};
%\draw [line width=0.3mm] (98,2.5) -- (98,2.6) node [above] {\scriptsize{$n = \tfrac{1}{2\eps}\mathrm{log}\tfrac{1}{\eps}$}};
%\begin{scope}
%\clip (0,-2.5) -- (0,2.5) -- (250,2.5) -- (250,-2.5) -- cycle;
%\draw [line width=0.3mm, dashed] (98,-2.5) -- (98,2.7);
%\draw[black] plot[only marks, mark=*,mark options={fill=white,xscale= 25*0.75*0.8,yscale=0.8333*0.75}] file {xoutvec_Fig2.txt};
%\draw[black] plot[only marks, mark=*,mark options={fill=black,xscale= 25*0.45*0.8,yscale=0.8333*0.45}] file {xvec_Fig2.txt};
%\end{scope}
%\end{tikzpicture}
\includegraphics[]{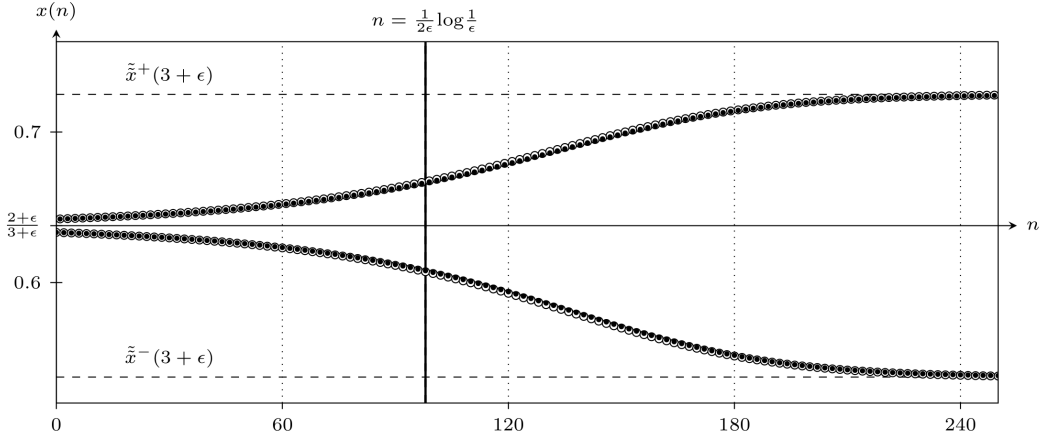}
\caption{Comparison between the numerical solution $x(n)$ and the two-term composite asymptotic expansion, which in this case is the late-time expansion given in (\ref{1.0 composite 2term}), for $\eps = 0.02$. The exact solution is indicated by filled circles, while the approximation is indicated by white circles. The point $n = \tfrac{1}{2\eps}\log(\tfrac{1}{2\eps})$ is indicated by a thick line. It is at this point that the early-time approximation breaks down; however, the composite approximation clearly provides an accurate approximation even once this point is passed.}\label{Static_late}
\end{figure*}

Converting back into the variables $x$ and $n$ of the original problem, the composite expansion based on \eqref{1.0 late-time} becomes
\begin{equation}
  \label{1.0 composite 2term}
  x(n) =  \frac{2}{3 } - \frac{(-1)^n\,\eps \, \e^{\eps n}}{3\,\sqrt{9 + \eps \, \e^{2\eps n} }} + \frac{18\,\eps - \eps^2 \,\e^{2\eps n}}{162 + 18 \, \eps \, \e^{2\eps n}} + \mathcal{O}(\eps^{\frac{3}{2}}).
\end{equation}
This composite approximation is compared with the numerical solution to the full problem in Figure \ref{Static_late} for the case where $\eps = 0.02$. While Figure \ref{F:StaticEarly} shows that the early time asymptotic approximations become invalid as $n$ grows large, here we see that the composite approximation remains valid over the entire domain.

\begin{figure}
\centering
%\begin{tikzpicture}
%[xscale=2,>=stealth,yscale=2]
%
%\draw[black] plot[only marks, mark=*,mark options={fill=white,xscale=0.25,yscale=0.25}] file {errfig_static1b.txt};
%\draw[black] plot[only marks, mark=*,mark options={fill=black,xscale=0.25,yscale=0.25}] file {errfig_static1.txt};
%
%\draw (-4.2,-2) -- (-0.8,-2) -- (-0.8,-5.5) -- (-4.2,-5.5) -- cycle;
%\draw[dotted] (-4.2,-3) -- (-0.8,-3);
%\draw[dotted] (-4.2,-4) -- (-0.8,-4);
%\draw[dotted] (-4.2,-5) -- (-0.8,-5);
%
%\draw (-4.2,-2) -- (-4.25,-2) node [left] {\scriptsize{-2}};
%\draw (-4.2,-3) -- (-4.25,-3) node [left] {\scriptsize{-3}};
%\draw (-4.2,-4) -- (-4.25,-4) node [left] {\scriptsize{-4}};
%\draw (-4.2,-5) -- (-4.25,-5) node [left] {\scriptsize{-5}};
%
%\draw (-4,-5.45) -- (-4,-5.55) node [below] {\scriptsize{-4}};
%\draw (-3,-5.45) -- (-3,-5.55) node [below] {\scriptsize{-3}};
%\draw (-2,-5.45) -- (-2,-5.55) node [below] {\scriptsize{-2}};
%\draw (-1,-5.45) -- (-1,-5.55) node [below] {\scriptsize{-1}};
%
%\node at (-2.5, -5.75) {\scriptsize{$\log_{10}(\eps)$}};
%\node at (-4.8, -3.5) {\scriptsize{$\log_{10}(\mathrm{error})$}};
%
%\draw (-2.75,-4.5+0.15) -- (-1.75,-3.5+0.15) -- (-1.75,-4.5+0.15) -- cycle;
%\node at (-2.25,-4.5+0.15) [below] {\scriptsize{$1$}};
%\node at (-1.75,-4+0.15) [right] {\scriptsize{$1$}};
%
%\filldraw (-3.8,-2.2) circle (0.02) ;
%\draw (-3.8,-2.5) circle (0.02);
%%\draw (-3.85,-2.8) -- (-3.75,-2.8);
%\node at (-3.7,-2.2) [right] {\scriptsize{Upper Branch Error}};
%\node at (-3.7,-2.5) [right] {\scriptsize{Lower Branch Error}};
%%\node at (-3.7,-2.8) [right] {\scriptsize{Slope = 1}};
%
%\end{tikzpicture}
\includegraphics[,width=0.65\columnwidth]{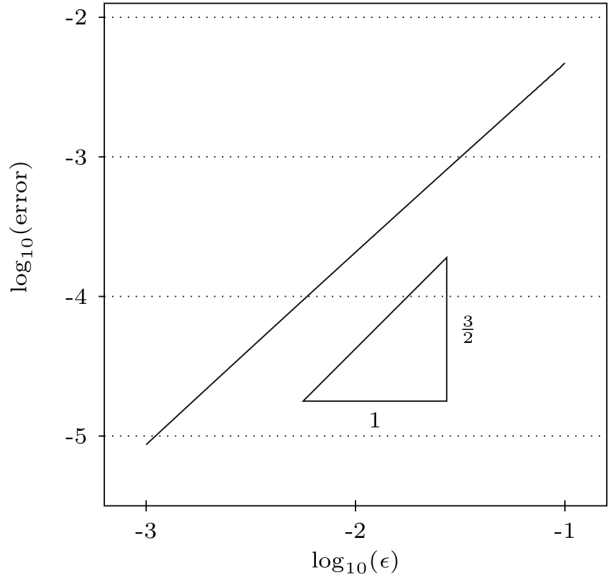}
\caption{Comparison of the approximation error of $x(n)$ as $\eps$ is varied. The approximation error is obtained by calculating $||x(n) - x_\mathrm{c}(n)||_{\inf}$, where $x_{\mathrm{c}}(n)$ is the two-term composite asymptotic approximation determined using (\ref{1.0 composite 2term}). For $\eps$ sufficiently small, the slope takes value $\frac{3}{2}$. Hence, the error is proportional to $\eps^{\frac{3}{2}}$, as predicted by the form of the composite expansion.}\label{1.0 staticerror}
\end{figure}

Figure \ref{1.0 staticerror} illustrates the approximation error of the composite expansion given in (\ref{1.0 composite 2term}) for a range of $\eps$ values, where the error is quantified as the infinity-norm of the difference between the exact solution $x(n)$ and the composite expansion. From Figure \ref{1.0 staticerror}, we see that the error is proportional to $\eps^{\frac{3}{2}}$ as $\eps \rightarrow 0$, which is consistent with the form of the asymptotic expansion (\ref{1.0 composite 2term}).

\section{Discrete logistic equation with $\bif = 3+\eps^2\,n$ and $x(0) = \frac{2}{3}$}
\label{S:Dynamic}

\subsection{Consequences of introducing a slowly varying bifurcation parameter}
\label{S:DynamicIntro}

In Section \ref{S:StaticAlgebraic}, we considered the behaviour of the discrete logistic equation when the bifurcation parameter was taken to be a constant, $ \bif = 3 + \eps$. In this section, we consider the dynamic problem described in \ref{S:IntroDiscLog}, where $\bif$ is taken to be a slowly varying function of $n$. Specifically, we will concentrate on a slow forward sweep through the bifurcation at $\bif = 3$, defining our small parameter, $\eps$, so that $\bif(n;\,\eps^2) = \bif_0 + \eps^2 \, n$, and hence the bifurcation parameter varies on the slow time scale associated with $\tau = \eps^2 \, n$. This definition of the time scale over which the bifurcation parameter evolves is for later convenience; as we will see, it will also be useful to define a slow time scale $t = \eps \, n$ associated with the solution trajectory.

With $\bif$ as given above, we find that the discrete logistic equation takes the form
\begin{equation}\label{0.0 Dynamic Temp}
x(n+1) = (\bif_0 +\eps^2 \, n) \, x(n) \, [1-x(n)].
\end{equation}
Before we apply our asymptotic methods to \eqref{0.0 Dynamic Temp}, it is useful to make a few observations about the behaviour of solutions to \eqref{0.0 Dynamic Temp} based on the analysis given in \cite{Baesens1991}.

Following Baesens \cite{Baesens1991}, we can define the adiabatic manifolds of \eqref{0.0 Dynamic Temp}, concentrating on the 1-periodic and 2-periodic manifolds. The 1-periodic adiabatic manifold, $\bar{x}(\bif;\eps^{2})$, is defined so that $\bar{x}(\bif)$ is invariant under \eqref{0.0 Dynamic Temp}, and so that $\bar{x}(\bif)$ is asymptotically close to $\xoneeqm(\bif)$, the 1-periodic equilibrium of the `static' discrete logistic equation defined in \eqref{1PeriodicSol}. Similarly, the two branches of the 2-periodic adiabatic manifold, $\bar{\bar{x}}^{\pm}(\bif;\eps^2)$, are defined so that each maps to the other under \eqref{0.0 Dynamic Temp}, and so that $\bar{\bar{x}}^{\pm}(\bif)$ is asymptotically close to $\xtwoeqm^{\pm}(\bif)$. Both adiabatic manifolds are obtained as asymptotic series in Section \ref{S:Adiabatic}.

When a parameter is changed dynamically through a bifurcation in which a steady state becomes unstable, it is usual for the bifurcation to appear delayed. This is because the solution trajectory has had an opportunity to become very close to the steady state before it becomes unstable, and it takes a reasonably long period of time to escape the unstable manifold. If the bifurcation parameter were then to be changed in the opposite direction, there would be a similar delay in the reverse bifurcation. This leads to the canard phenomenon, first described in \cite{Benoit1981}, with the particular case of discrete canards discussed in \cite{Fruchard1988,Fruchard2003}. Baesens \cite{Baesens1991} discussed the canard problem for the discrete logistic equation, using renormalisation methods to analyse both forward sweep and backward sweep. In this paper, we concentrate solely on the more challenging problem of forward sweep, but the methods described here could also be used to obtain trajectories associated with backward sweep.

If we were to pick $\bif_0 < 3$ so that $n = \ord(\eps^{-2})$ when $\bif = 3$, we would find that $x(n)$ is exponentially close to the 1-periodic adiabatic manifold when $\bif = 3$. As discussed in \cite{Baesens1991} and noted in Section \ref{S:StaticAlgebraic}, this has the consequence that it would then take a $\tau = \ord(1)$ amount of time for $x(n)$ to escape the neighbourhood of the 1-periodic adiabatic manifold after it has become unstable. %Indeed, if we were to pick $\bif_0 < 2$, we would not see any of the bifurcations associated with the period doubling cascade to chaos in our solution to \eqref{0.0 Dynamic Temp}, since $\bif$ would exceed 4 (meaning that trajectories are unbounded) before $x(n)$ escapes the neighbourhood of the 1-periodic adiabatic manifold.

However, this severe delay to the bifurcation is very delicate, and can only be observed in the absence of any noise or rounding error in the numerical calculation of the solution to \eqref{0.0 Dynamic Temp}. Even a very small amount of noise will mean that $x$ is not exponentially close to $\bar{x}(\bif)$ when $\bif \geq 3$, and the calculated solution will rapidly escape the unstable 1-periodic adiabatic manifold. A significant part of the analysis in \cite{Baesens1991} is devoted to understanding the effect of noise on the solution to the discrete logistic equation with a slowly varying bifurcation parameter.

In order to avoid the problems associated with the effects of noise (and thereby allow clear numerical validation for our asymptotic results), we will concentrate on the case where $\bif_0 = 3$, and $x(0)$ is algebraically close to the critically stable 1-periodic adiabatic manifold. Specifically, we concentrate on the case where $x(0) = \frac{2}{3}$, and hence $x(0) - \bar{x}(0) = \ord(\eps^{2})$. 

\begin{figure*}
\centering
\includegraphics[]{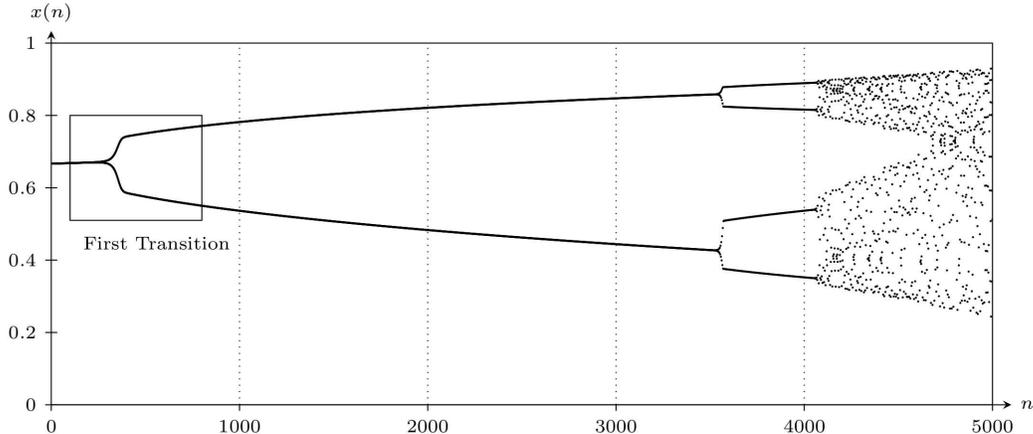}
\caption{Numerical solution to \eqref{0.0 Dynamic}, the discrete logistic equation with slowly varying bifurcation parameter, when $\eps = 0.012$.}\label{Dynamic_exact}
\end{figure*}

Thus, the focus of this section will be on analysing the non-autonomous difference equation specified by
%\begin{subequations}
\begin{equation}\label{0.0 Dynamic}
x(n+1) = (3 +\eps^2 \, n) \, x(n) \, [1-x(n)],
\end{equation}
with $x(0) = \tfrac{2}{3}$.
%subject to the initial condition $x(0) = \tfrac{2}{3}$.
%\begin{equation}
% \label{Dynamic IC}
% x(0) = \frac{2}{3}.
%\end{equation}
%\label{DynamicUnscaled}
%\end{subequations}
A numerical solution to this problem with $\eps = 0.012$ is illustrated in Figure \ref{Dynamic_exact}. Note that all of the bifurcations are delayed, and there are rapid transitions between the adiabatic manifolds. Note also that the numerical solution cannot be treated as reliable once we reach the transition from the 2-periodic manifold to the 4-periodic manifold. %When $\lambda$ reaches the bifurcation value of $1 + \sqrt{6}$, $x(n)$ should be exponentially close to the critically stable 2-periodic adiabatic manifold. Hence, tiny (but only algebraically small) amounts of numerical noise could make $x(n)$ further away from the 2-periodic manifold at this critical point, potentially leading to significant changes in when the solution leaves the 2-periodic manifold and approaches the 4-periodic manifold. 

\subsection{Adiabatic manifolds and early time scaling}
\label{S:Adiabatic}

In Section \ref{S:StaticAlgebraicIntro}, we simply rescaled $x$ around $\tfrac{2}{3}$ in order to obtain a weakly nonlinear problem for the escape from the 1-periodic equilibrium with a static bifurcation parameter. In our analysis of the discrete logistic equation with a slowly varying bifurcation parameter, we wish to instead rescale about the 1-periodic adiabatic manifold, $\bar{x}(3 + \eps^2 \, n)$. This means that we need to obtain an asymptotic expression for the 1-periodic adiabatic manifold before we seek asymptotic solutions to \eqref{0.0 Dynamic}.

To find the 1-periodic adiabatic manifold, we start from \eqref{0.0 Dynamic Temp}, which we can rewrite in the form
\begin{equation}\label{FindingAdiabatic1 Temp}
x(n+1) = \lambda(\tau) \, x(n) \, [1-x(n)],
\end{equation}
where $\tau = \eps^2 \, n$ and $\diff{\lambda}{\tau} = 1$.

Since the 1-periodic adiabatic manifold is dependent only on $\lambda$, and $\diff{\lambda}{\tau} = 1$, we can reparameterise \eqref{FindingAdiabatic1 Temp} in terms of $\lambda$, so that it takes the form
\begin{equation}
 \bar{x}(\lambda + \eps^2) = \lambda \, \bar{x}(\lambda) \, [1 - \bar{x}(\lambda)].
\end{equation}
Using Taylor series to expand the left hand side around $\eps = 0$ in powers of $\eps^2$, and solving the resulting system of algebraic equations leads to
%\begin{equation}\label{2.1 stable1}
%\sum_{j=0}^{\inf} \frac{\eps^{2j}}{j!}\diff{^j \bar{x}(\bif)}{\bif^j} =  \bif \, \bar{x}(\bif) \, [1-\bar{x}(\bif)].
%\end{equation}
%We now propose an asymptotic series expansion for $\bar{x}(\bif)$ in even powers of $\eps$, and exploit the fact that $\bar{x}(\bif) = \xoneeqm(\bif) + \Ord(\eps^2)$. Collecting like powers of $\eps$ in \eqref{2.1 stable1}, we then find that the 1-periodic adiabatic manifold is
\begin{multline}\label{1Adiabatic}
 \bar{x}(\bif) 
 = \frac{\bif-1}{\bif} 
 - \eps^2 \, \frac{1}{(\bif - 1)\, \bif^2} 
 + \eps^4 \, \frac{\bif^2 - 5 \, \bif + 2}{(\bif-1)^3 \, \bif^3} \\
 + \eps^6 \, \frac{-\bif^4 + 14 \, \bif^3 - 44 \, \bif^2 + 27 \, \bif - 6}{(\bif-1)^5 \, \bif^4}
 + \Ord(\eps^8).
\end{multline}
Further terms can be obtained by matching at higher powers of $\eps$.

%The 2-periodic adiabatic manifold can similarly be obtained by considering the doubled map, which takes the form
%\begin{multline}\label{FindingAdiabatic2 Temp}
%x(n+2) = \bif(\tau) \, \bif(\tau+ \eps^2) \, x(n) \, \big[1-x(n)\big] \\
% \times \big(1 - \bif(\tau) \, x(n) \, \big[1 - x(n)\big] \big),
%\end{multline}
%where again $\diff{\bif}{\tau} = 1$. Following the same procedure of reparameterising in terms of $\bif$ and expanding the series with $\bar{\bar{x}}^\pm(\bif) = \xtwoeqm^\pm(\bif) + \Ord(\eps^2)$, we find that
Starting from the doubled map and repeating the procedure yields
\begin{equation} \label{2Adiabatic}
 \bar{\bar{x}}^{\pm}(\bif)
 = \frac{\bif + 1 \pm \sqrt{\bif^2 - 2 \bif -3}}{2 \bif} 
 + \eps^2 \bigg[\frac{\bif + 3}{2 \, (\bif-3)\,\bif^2\,(\bif+1)}
 \pm \frac{\bif^2 - 4 \, \bif - 9}{2 \, (\bif-3)^{\frac{3}{2}}\,\bif^2\,(\bif+1)^{\frac{3}{2}}}
 \bigg] 
 + \mathcal{O}(\eps^4).
\end{equation}
Other adiabatic manifolds could be obtained using the same method.

%The 2-periodic adiabatic manifold will be considered further in Section \ref{S:DynamicFailure}.

% \subsection{Early time scaling}
% \label{S:DynamicScaling}

Having obtained the 1-periodic adiabatic manifold in \eqref{1Adiabatic}, we can rescale $x$ in \eqref{0.0 Dynamic} to obtain a weakly nonlinear problem by introducing $X(n)$ so that
\begin{equation}
 \label{DynamicScaling Temp}
 x(n) = \bar{x}(3 + \eps^2 \, n) + \eps^2 \, X(n).
\end{equation}
However, while this would lead to a rescaled version of \eqref{0.0 Dynamic} that is appropriate when $n = \Ord(1)$, we wish to use the method of multiple scales to find $X(n)$ for $n \gg 1$, which means that we need to introduce a slow time scale, $t$, alongside $n$. By comparing the variation of the bifurcation parameter $\lambda(n;\eps) = 3 + \eps^2 n$ with the form of the Taylor series expansion for $X(n+1)$, we find that the appropriate slow time scale is again given by $t = \eps n$.

%
%Since \eqref{0.0 Dynamic} and \eqref{DynamicScaling Temp} are both dependent on $n$, this has the potential to affect the form of the rescaled asymptotic expansions. 

%Introducing a general slow time scale, $t = \eps^p \, n$ for some value of $p$, we firstly note that $X(n+1) \equiv X(n+1,t+\eps^p)$ will be expressed as a Taylor series in powers of $\eps^p$. At the same time, the fact that the bifurcation parameter is varying on the slow time scale means that it makes most sense to write $\bif(n;\,\eps) = 3 + \eps^2 \, n$ in the form $\bif = 3 + \eps^{2-p} \, t$. In order for the derivatives from the Taylor series expansion and the terms associated with the slowly varying bifurcation parameter to appear at equivalent orders in $\eps$, we therefore choose $p = 2-p$, and so $p = 1$.

Thus, our early time scaling for the method of multiple scales takes the form
\begin{equation}
 \label{DynamicScaling}
 x(n) = \bar{x}(3 + \eps \, t) + \eps^2 \, X(n,\,t).
\end{equation}
Using \eqref{1Adiabatic} and introducing the scaling from \eqref{DynamicScaling} to \eqref{0.0 Dynamic} yields
%\begin{subequations}
%\label{DynamicScaled}
\begin{equation}
 X(n+1,\,t+\eps) + X(n,\,t) 
 = - \eps \, t \, X(n,\,t) 
 + \eps^2 \left[\tfrac{1}{3} \, X(n,\,t) - 3 \, X(n,\,t)^2  \right]
 + \Ord(\eps^3),
 \label{DynamicScaled Main}
\end{equation}
which must be solved subject to
\begin{equation}
 X(0,0) = \tfrac{1}{18} + \tfrac{1}{54}\,\eps^2 + \Ord(\eps^4). \label{DynamicScaled ICs}
\end{equation}
%\end{subequations}
More terms can easily be obtained by using higher orders in the asymptotic expansion of $\bar{x}(\bif)$ given in \eqref{1Adiabatic}, noting throughout that we are using $\bif = 3 + \eps \, t$.

\subsection{Early time solution}
\label{S:DynamicEarly}

We now express $X(n,t)$ as a power series in $\eps$ and expand $X(n+1,t+\eps)$ using Taylor series as previously. Substituting into \eqref{DynamicScaled Main} and collecting like orders of $\eps$, we obtain
\begin{align}
\mathcal{O}(1)&: %\notag \\
 & X_0(n+1,t) + X_0(n,t) &= 0, \label{2.2 one} \\
\mathcal{O}(\eps)&: %\notag \\
 & X_1(n+1,t) + X_1(n,t) &=  - t \, X_0(n,t) - \diff{X_0(n,t)}{t}, \label{2.2 eps}
%\mathcal{O}(\eps^2)&: %\notag \\
% & X_2(n+1,t) + X_2(n,t) &= - t \, X_1(n,t)  \notag \\
%& &  + \frac{1}{3} \, X_0(n,t) - 3 \, &X_0(n,t)^2 - \diff{X_1(n,t)}{t} \notag \\
%& & & \, - \frac{1}{2} \diff{^2 X_0(n,t)}{t^2},\label{2.2 eps2}
\end{align}
which must be solved subject to initial conditions based on \eqref{DynamicScaled ICs}.

Solving \eqref{2.2 one} yields $X_0(n,t) = (-1)^n A(t)$. Continuing to the next order, we obtain the secularity condition $A'(t) - t \, A(t) = 0$; combined with the initial condition $X_0(0,0) = A(0) = \frac{1}{18}$, this yields
%
%
%Solving this equation with the appropriate initial condition gives $X_0(n,t) = 0$. The expression in (\ref{2.1 eps1}) now simplifies to give
%\begin{equation}
%X_1(n+1,t) + X_1(n,t) = 0,
%\end{equation}
%with the condition $X_1(0,0) = 1/18$. Solving this gives $X_1(n,t) = (-1)^n B(t)$, with $B(0) = \tfrac{1}{18}$. At the subsequent order (\ref{2.2 eps2}), we obtain the secularity condition as $B'(t) - t B(t) = 0$. Solving this with the appropriate boundary condition gives
\begin{equation}\label{Dynamic EarlyLeading}
X_0(n,t) = \tfrac{1}{18}(-1)^n\e^{\frac{t^2}{2}}.
\end{equation}
This process can be continued to higher orders, and we obtain
\begin{equation}\label{2.0 threeterm}
X(n,t) = \tfrac{1}{18}(-1)^n\e^{\frac{t^2}{2}} 
- \tfrac{\eps}{108} \, (-1)^n \, \left({t^3} + {5t}\right)\e^{\frac{t^2}{2}} 
+ \mathcal{O}(\eps^2).
\end{equation}

\begin{figure}
\centering
%\begin{tikzpicture}
%[xscale=0.04*1*0.9*0.5,>=stealth,yscale=1.2*0.8*0.9*0.5]
%
%\draw[->] (0,0) -- (315,0) node[right]{\scriptsize{$n$}};
%\draw[->] (0,6.2) -- (0,6.6) node[above]{\scriptsize{$x_n$}};
%\draw (4,6) -- (-4,6) node[left]{\scriptsize{$6$}};
%\draw (4,3) -- (-4,3) node[left]{\scriptsize{$3$}};
%\draw (4,0) -- (-4,0) node[left]{\scriptsize{$0$}};
%\draw (4,-3) -- (-4,-3) node[left]{\scriptsize{$-3$}};
%\draw (4,-6) -- (-4,-6) node[left]{\scriptsize{$-6$}};
%\draw [dotted] (100,-6.2) -- (100,6.2);
%\draw [dotted] (200,-6.2) -- (200,6.2);
%\draw [dotted] (300,-6.2) -- (300,6.2);
%%\draw [dotted] (160,-2.5) -- (160,2.5);
%\draw (100,-6.2) -- (100,-6.4) node[below] {\scriptsize{$40$}};
%\draw (200,-6.2) -- (200,-6.4) node[below] {\scriptsize{$80$}};
%\draw (300,-6.2) -- (300,-6.4) node[below] {\scriptsize{$120$}};
%%\draw (160,-2.5) -- (160,-2.6) node[below] {\scriptsize{$160$}};
%\draw (0,-6.2) -- (0,-6.4) node[below] {\scriptsize{$0$}};
%
%\begin{scope}
%\clip (0,-6.2) -- (0,6.2) -- (301,6.2) -- (301,-6.2) -- cycle;
%
%%\draw[black] plot[only marks, mark=*,mark options={fill=white,xscale= 25*0.5*0.8,yscale=0.8333*0.5}] file {xinn1vec_Fig1.txt};
%\draw[black] plot[only marks, mark=*,mark options={fill=white,xscale= 25*0.5*0.8*2.5,yscale=0.8333*0.5*2.5}] file {xinn2vec_Fig3.txt};
%\draw[black] plot[only marks, mark=*,mark options={fill=black,xscale= 25*0.5*0.8*1.5,yscale=0.8333*0.5*1.5}] file {xvec_Fig3.txt};
%
%\end{scope}
%
%\draw (0,-6.2) -- (0,6.2) -- (301,6.2) -- (301,-6.2) -- cycle;
%
%
%\end{tikzpicture}
\includegraphics[]{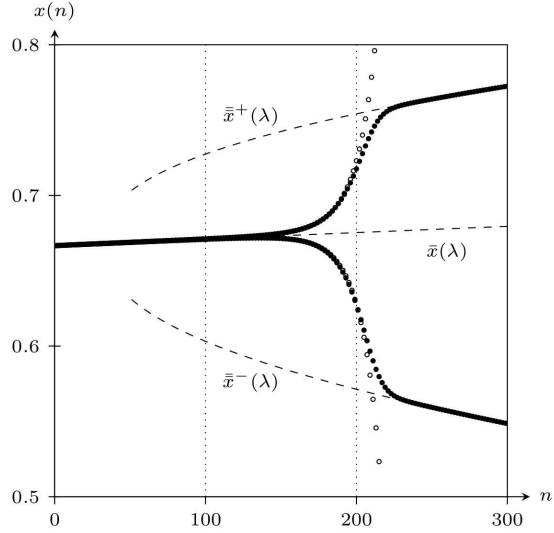}
\caption{Comparison between the numerical solution $x(n)$ the early-time approximation (\ref{2.0 threeterm}), for $\eps = 0.02$. The exact solution is represented by black circles, while the asymptotic approximation is represented by white circles. As in the static case, the early-time solution provides an accurate approximation until the transition region is reached, at which point it diverges from the exact value.}\label{Sweeping_early}\end{figure}

This approximation is illustrated in Figure \ref{Sweeping_early} for $\eps = 0.02$. We see again that, although this approximation is accurate for early times, it loses accuracy when $t$ becomes sufficiently large, and hence does not capture the transition to the 2-periodic adiabatic manifold. As before, however, this can be captured using the method of matched asymptotic expansions.

\subsection{Failure of the early time expansion}
\label{S:DynamicFailure}

The first step in determining the late time asymptotic solution to \eqref{0.0 Dynamic} is to determine when the early time expansion fails, as in Section \ref{S:StaticAlgebraicFailure}. 

Taking the doubled map and assuming that $x$ can be expressed as a function of $t = \eps \, n$ only, we now introduce the early time ansatz
%so that $x(n) \equiv \hat{\chi}(t)$, we see that the doubled map takes the form
%\begin{multline}
% \hat{\chi}(t + 2 \, \eps) = (3 + \eps \, t) \, (3 + \eps \, t + \eps^2 ) \, \hat{\chi}(t) \, \big[1-\hat{\chi}(t)\big] \\
% \times \big(1 - \big[3 + \eps \, t \big] \, \hat{\chi}(t) \, \big[1 - \hat{\chi}(t)\big] \big).
% \label{DynamicChi Temp}
%\end{multline}
\begin{equation}
x(n) \equiv \bar{x}(3+\eps \, t) + \eps^2 \, \chi(t),
\end{equation} 
where $\bar{x}(\bif)$ is obtained from \eqref{1Adiabatic}. We will substitute this ansatz into the doubled map and seek a new dominant balance that is consistent with the early time solution, $\chi(t) = \ord[\exp (\frac{t^2}{2})]$. Noting that we are only concerned with the largest terms at each order of $\eps$ (\ie{} the terms with the largest power of $\chi$, since $\chi$ is exponentially growing in $t$), this yields
\begin{align}
  \nonumber \underset{\textrm{\ding{172}}}{2 \, \diff{\chi}{t}} + \sum_{j=2}^{\inf}\frac{\eps^{j-1} \, 2^j}{j!} \, \diff{^j\chi}{t^j} 
  = & \, \underset{\textrm{\ding{173}}}{2 \, t  \, \chi(t)} 
  + \eps \, \big[t^2 \, \chi(t) + \ldots\big] 
 - \eps^2 \, \big[3 \,  t \, \chi(t)^2  + \ldots\big] 
 - \eps^3 \, \big[\underset{\textrm{\ding{174}}}{18 \, \chi(t)^3} + \ldots\big] \\ 
 & - \eps^4 \, \big[30 \, t \, \chi(t)^3 + \ldots\big]
 - \eps^5 \, \big[27 \, \chi(t)^4 + \ldots \big] 
 + \Ord[\eps^6 \, \chi(t)^4 \, t^\nu],
 \label{DynamicChi}
\end{align}
where $\nu$ is finite, where the omitted terms are at most $\Ord[\eps^6 \, \chi(t)^4 \, t^\nu]$, and will not contribute to the change in dominant balance.

As in Section \ref{S:StaticAlgebraicFailure}, we propose rescalings $\chi = \delta \, \xi$ and $t = K_0 + K_1 \, s$, where $\xi$ and $s$ are $\mathcal{O}(1)$ variables when the new dominant balance is attained, and $\delta \gg 1$, $K_0 \gg 1$ and $K_0 \gg K_1$.  We now seek a dominant balance in \eqref{DynamicChi} that is consistent with the leading-order behaviour of $\chi(t)$ based on \eqref{Dynamic EarlyLeading}. 

For consistency with \eqref{Dynamic EarlyLeading}, we require
\begin{equation}
\delta \, \xi = \ord\left[\exp \left(\tfrac{1}{2}{{K_0}^2} + K_0 \, K_1 \, s + \tfrac{1}{2}{{K_1}^2 \, s^2}\right) \right].
\end{equation}
Since this must hold when $\xi$ and $s$ are $\mathcal{O}(1)$, it is appropriate to choose $K_1 = K_0^{-1}$ and $\delta = \exp[ \frac{{K_0}^2}{2}]$. 

Writing $K_0$ as $K$, we make the substitution $t = K + K^{-1} \, s$ in \eqref{DynamicChi}, and note that $\diff{}{t} \equiv K \, \diff{}{s}$. As a result of this, we see that the terms labelled \ding{172} and \ding{173} will always be the same size, $\Ord(\delta \, K)$, and will always participate in the dominant balance.

As $t$ and $\chi$ increase, the first term to become the same size as \ding{172} and \ding{173} is the term labelled \ding{174}. Thus, the new dominant balance occurs when $\eps^3 \, \delta^3 = \ord[K \, \delta]$, and we pick $K$ to be a solution of the transcendental equation
\begin{equation}
\eps^3 K^{-1}\e^{K^2} = 1. \label{KTranscendental}
\end{equation}
Noting that $K$ must be large as $\eps \rightarrow 0$, this equation may be solved in terms of the $-1$ branch of the Lambert W function (see, for example, \cite{LambertW}). For our purposes, however, it is more useful to define an iterative procedure for the solution by noting that
\begin{equation}
K = \sqrt{3 \log \left(\tfrac{1}{\eps}\right) + \log K},
\end{equation}
and hence $K$ can be expressed as an asymptotic series as $\eps \rightarrow 0$, with the leading order behaviour of $K$ being given by $K \sim \sqrt{3} \, \log(\frac{1}{\eps})^\frac{1}{2}$.

Since $\delta = \exp[ \frac{{K_0}^2}{2}]$, we note that \eqref{KTranscendental} gives
\begin{equation}
 \delta = K^\frac{1}{2} \, \eps^{-\frac{3}{2}}.
\end{equation} 
In the analysis of the late time solution that follows, we will treat $K$ and $\eps$ as independent parameters, although they are linked by \eqref{KTranscendental}. % Since $K$ is only logarithmically large in $\eps$, we will treat $K$ as if it were an $\Ord(1)$ parameter for the purposes of asymptotic analysis.

\subsection{Late time solution and matched asymptotic expansions}

Having found that $K$ is logarithmically large as $\eps \rightarrow 0$ and having obtained an expression for $\delta$ in terms of $\eps$ and $K$, we can rescale our original equation, \eqref{0.0 Dynamic}, by introducing late time variables $\xi$, $m$ and $s$, so that
%\begin{subequations}
%\label{DynamicLateRescaling}
\begin{align}
 x(n) &= \bar{x}(3 + \eps \, K + \eps \, K^{-1} \, s) + \eps^\frac{1}{2} \, K^\frac{1}{2} \, \xi(m,\,s), \label{DynLate1} \\
 t &= K + K^{-1} \, s \label{DynLate2} \\
 n &= \tfrac{K}{\eps} - \gamma + m, \label{DynLate3}
\end{align}
%\end{subequations}
where $\gamma$ is chosen so that $0 \leq \gamma < 2$, and $m \equiv n \pmod 2$ for later convenience.

Substituting \eqref{DynLate1}, \eqref{DynLate2}, and \eqref{DynLate3} into \eqref{0.0 Dynamic} and using \eqref{1Adiabatic}, we find as the rescaled form of \eqref{0.0 Dynamic} for late times
\begin{multline}
 \xi(m+1,\,s) + \xi(m,\,s)
 = - \eps^\frac{1}{2} \, K^{\frac{1}{2}} \, \xi(m,s)^2 \\
 - \eps \, K \,  \left(\left[1 + \frac{s}{K^2}\right] \, \xi(m,s) 
 + \diff{\xi(m,s)}{s} \right)
 + \Ord(\eps^\frac{3}{2} \, K^{\frac{3}{2}}).
 \label{DynamicLate}
\end{multline} 

We now require a series expansion for $\xi(m,\,s)$.
Since $K\rightarrow \infty$ as $\eps \rightarrow 0$, one way of proceeding would be to express each $\xi(m,\,s)$ as a double asymptotic series in $\eps$ and $K$. Inspection of \eqref{DynamicLate} indicates that an appropriate ansatz would be
\begin{equation}
 \xi(m,\,s) \sim \sum_{k=0}^\infty \sum_{l=0}^\infty \eps^{\frac{k}{2}} \, K^{\frac{k}{2}-2l} \, \xi_{k,l}(m,\,s),
\end{equation}
however, the fact that $K$ is only logarithmically large means that it is more convenient to treat $K$ as though it were an $\Ord(1)$ parameter, and write this expression as a single series in powers of $\eps$.

%However, the fact that $K$ is only logarithmically large means that it is more convenient to treat $K$ as though it were an $\Ord(1)$ parameter, and express $\xi(m,\,s)$ as a single asymptotic series. This is advantageous because it would otherwise be necessary to obtain $\xi(m,\,s)$ at all algebraic orders of $K$ before it is possible to move to the next order in $\eps$. Moreover, expanding in the logarithmically large parameter $K$ would make it more difficult to use Van Dyke's rule for the method of matched asymptotic expansions. Hence, we instead express $\xi(m,\,s)$ as a series of the form
%\begin{equation}
% \xi(m,\,s) \sim \sum_{k=0}^{\infty} \eps^{\frac{k}{2}} \, K^{\frac{k}{2}} \, \xi_k(m,\,s). \label{DynamicXiSeries}
%\end{equation}

Substituting into \eqref{DynamicLate} and collecting $\Ord(1)$ terms, we find that $\xi_0(m,\,s) = (-1)^m \, P(s)$ for some smooth function $P(s)$ to be determined. Continuing to $\Ord(\eps^\frac{1}{2})$ gives information about $\xi_1(m,\,s)$, but does not determine $\xi_0(m,\,s)$. Then, at $\Ord(\eps)$, we obtain the secularity condition
\begin{equation} \label{DynamicLateSlow}
P'(s) =  P(s) + \frac{s \, P(s)}{K^2} - 9 \, P^3(s),
\end{equation}
which yields
\begin{equation}\label{2.2 B0}
P(s) = \frac{\e^{s + \frac{s^2}{2 K^2}}}{\sqrt{\kappa_p + 18 \,K \, \exp\left[2\,s + \frac{s^2}{K^2}\right] \, F\left[K + \frac{s}{K} \right]}},
\end{equation}
where $\kappa_p$ is an arbitrary constant, and $F(x)$ is Dawson's integral, defined by
\begin{equation}
 F(x) = \e^{-x^2} \int_0^x \e^{y^2} \, \d y.
\end{equation}
%Note that \eqref{2.2 B0} can easily be expanded as a regular asymptotic series in the small parameter $K^{-2}$, yielding the series expansion that would be obtained using the alternative method outlined above.

In order to determine $\kappa_p$ using Van Dyke's rule, we express \eqref{2.2 B0} in terms of the early time variables and take the leading order term, which yields the result that
\begin{equation}
 K^\frac{1}{2} \, \eps^{-\frac{3}{2}} \, \xi_0 \sim (-1)^n \, \e^\frac{t^2}{2}{\kappa_p^{-\frac{1}{2}}}.
\end{equation}
Since this must match with the early time solution \eqref{Dynamic EarlyLeading} expressed in late time variables, we find that $\kappa_p = 324$, and hence the leading order late time solution is given by
\begin{equation}\label{2.2 B02}
\bar{\xi}_0(m,s) =  \frac{(-1)^m \, \e^{s + \frac{s^2}{2 \, K^2}}}{\sqrt{324 + 18 \,K \, \exp\left[2\,s + \frac{s^2}{K^2}\right] \, F\left[K + \frac{s}{K} \right]}}.
\end{equation}

As previously, the matching term is identical to the early time solution, and so the uniformly valid composite approximation is identical to the late time solution. Rewriting \eqref{2.2 B02} in terms of the original variables and using the first two terms from the series expansion of $\bar{x}(\bif)$ from \eqref{1Adiabatic}, we find that the composite approximation is given by
\begin{equation}
\label{Dynamic Composite}
 x(n) \sim \frac{2 + \eps^2 \, n}{3 + \eps^2 \, n} 
 - \eps^2 \, \frac{1}{(2+\eps^2 \, n)(3+\eps^2 \, n)^2} 
 + \eps^2 \, \frac{(-1)^n \, \e^{\frac{\eps^2 n^2}{2}}}{\sqrt{324 + 9 \, \pi^\frac{1}{2} \, \eps^{\frac{3}{2}} \erfi[\eps \, n]}},
\end{equation}
where $\erfi(x)$ is the imaginary error function, defined by
\begin{equation}
 \erfi(x) = \frac{\erf(\mathrm{i}\,x)}{\mathrm{i}} = \frac{2}{\sqrt{\pi}} \int_0^{x} \e^{y^2} \, \d y.
\end{equation} 
We note that \eqref{Dynamic Composite} will be accurate up to $\Ord(\eps^3)$ in the early time region, and accurate up to $\Ord[\eps \, K]$ (or, equivalently, $\Ord[\eps \log(\frac{1}{\eps})^\frac{1}{2}]$) in the late time region. 

\begin{figure*}
\centering

\includegraphics[]{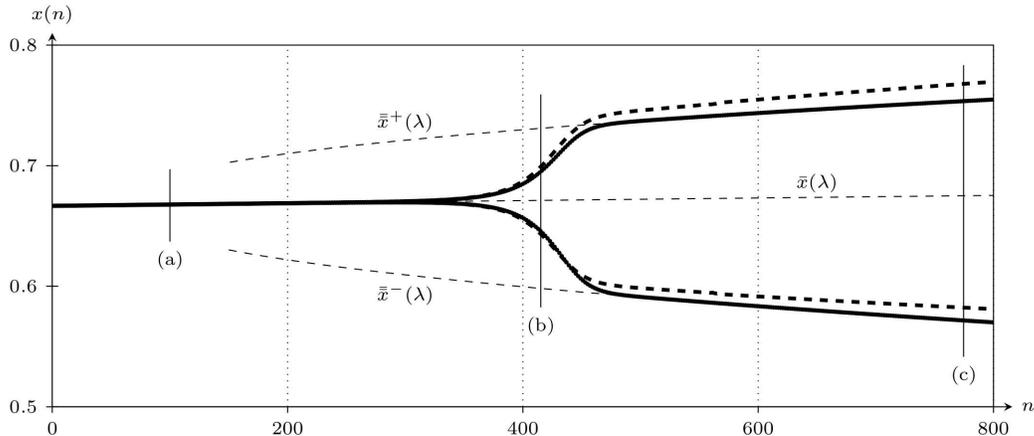}
\caption{Comparison between the numerical solution $x(n)$ for $\eps = 0.01$, denoted by dense black circles, and the leading-order composite expansion given in (\ref{Dynamic Composite}), with the two branches of the asymptotic solution illustrated by thick dashed lines. We expect the error to be $\mathcal{O}(\eps^3)$ in the early-time region, and $\mathcal{O}[\eps \log(\tfrac{1}{\eps})^{\frac{1}{2}}]$ in the late-time region. The points (a), (b) and (c) indicate representative examples of the early-time, transition region, and late-time behaviour, and correspond to Figures \ref{fig dynamicerror} (a), (b) and (c) respectively.}\label{Sweeping_comp}
\end{figure*}

In Figure \ref{Sweeping_comp}, we compare the composite approximation from \eqref{Dynamic Composite} with the numerical solution to the full problem when $\eps = 0.01$. This shows very good agreement on the position and width of the rapid transition between the early time behaviour and the late time behaviour, but we also see that \eqref{Dynamic Composite} is much more accurate at early times than at late times. 

\begin{figure*}
\centering
%\subfloat[$n = \left\lfloor\tfrac{1}{\eps}\right\rfloor$]{
\includegraphics[width=0.999\textwidth]{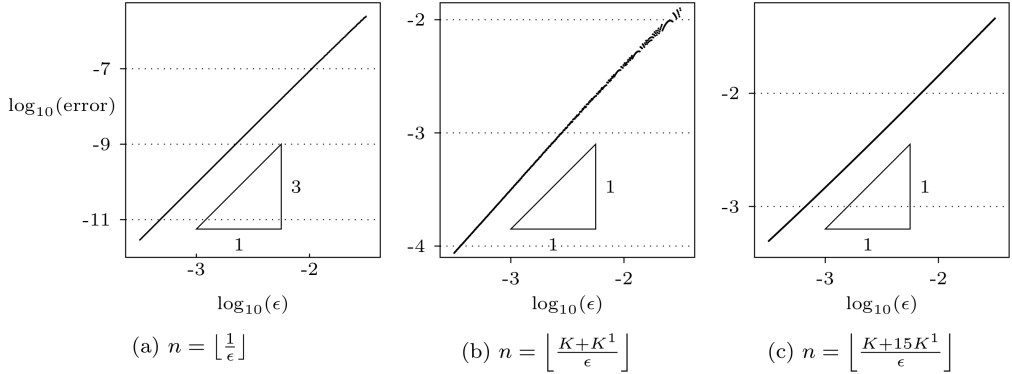}
%}
%\subfloat[$n \approx \tfrac{K+K^{1}}{\eps}$]{
%\includegraphics[width=0.3\textwidth]{Figures_FIG12b.eps}
%}
%\subfloat[$n \approx \tfrac{K+15 K^{1}}{\eps}$]{

%\includegraphics[width=0.3\textwidth]{Figures_FIG12c.eps}
%}
\caption{Comparison of the approximation error of $x(n)$ as $\eps$ is varied in (a) the early-time solution, (b) the transition region, and (c) the late-time solution. The approximation error is given by calculating the difference between the exact solution and the asymptotic approximation at particular values of $n$ in each region. In (a), the error is calculated at $n = \lfloor \frac{1}{\eps} \rfloor$, which corresponds to $t = \mathcal{O}(1)$. In (b), the error is calculated at $n \approx (K+K^{-1}) \, \eps^{-1}$, which corresponds to $s = \mathcal{O}(1)$. Finally, in (c), the error is calculated at $n \approx (K + 15 K^{-1}) \, \eps^{-1}$, which corresponds to $s$ increasing out of the transition region and into the late-time behaviour. In each of these cases, the $n$ used for calculation is taken to be the nearest value of $n$ that lies on the solution branch with the greater approximation error.  These three points are illustrated on the composite approximation shown in Figure \ref{Sweeping_comp}.  The unusual behaviour in (b) as $\eps$ becomes less small is a consequence of the restricted choice of $n$ due to the discrete nature of the fast scale. } % DISCUSS IN TEXT We find that the slope in (a) takes the value three, showing that the error is proportional to $\eps^3$, which is consistent with the early-time solution given in [REF]. In (b) and (c), the slope takes unit value, indicating that the error is proportional to $\eps$. This demonstrates that the composite asymptotic expansion provides an accurate approximation to both the early-time and late-time behaviour, as well as capturing the behaviour within the transition region.}
\label{fig dynamicerror}
\end{figure*}

Because the accuracy of our approximation varies considerably with $n$, we perform an error analysis by considering the difference between the numerical solution and the composite approximation at three illustrative points, marked (a), (b) and (c) on Figure \ref{Sweeping_comp}. Keeping a fixed value of $t$ in the case of (a) and a fixed value of $s$ in the case of (b) and (c), we consider the size of the error as $\eps \rightarrow 0$. As shown in Figure \ref{fig dynamicerror}, we see that the error appears to be $\Ord(\eps^3)$ for fixed $t$ and $\Ord(\eps)$ for fixed $s$. This is consistent with our expectations, since the logarithmic term in the late time error will be very difficult to observe.

%The method of matched asymptotic expansions is most commonly associated with singular perturbation problems in differential equations where there is a small parameter multiplying the highest derivative. In the matched asymptotic expansions that we consider in this paper, we will not see any analogous small parameter multiplying the `highest difference'. This is not a problem, however, since the method of matched asymptotic expansions is well known to be more general than this; as presented by both Hinch \cite{HinchPert} and Nayfeh \cite{NayfehPert}, matched asymptotic expansions can be used whenever there is a change in the dominant balance from one scaling regime to another. One unusual feature of the 

\section{Discussion and conclusions}
\label{S:Discussion}

\subsection{Strategies for the asymptotic analysis of difference equations}
\label{S:Strategy}

In this paper, we have demonstrated that a combination of the method of multiple scales and the method of matched asymptotic expansions can be used to obtain uniformly valid asymptotic solutions to the discrete logistic equation in the neighbourhood of a period doubling bifurcation. We have developed a method for obtaining the late time rescaling by eliminating the effects of the discrete fast time scale for both autonomous and non-autonomous difference equations. We have concentrated on examples that arise from particular individual difference equations, but the same approach could also be used for systems of difference equations.

In the present work, we have not made any attempt to define a class of singularly perturbed difference equations where it is possible to prove that the combination of multiple scales and matched asymptotic expansions will lead to a uniformly valid composite asymptotic solution. It seems highly unlikely that such a proof would be straight forward to obtain or particularly enlightening; as discussed in detail by Eckhaus \cite{Eckhaus1979,Eckhaus1994} and noted elsewhere (see, for example, \cite{Fraenkel1969-I,Lagerstrom-MAE}), it is notoriously difficult to develop a rigorous and general theory to justify the method of matched asymptotic expansions (and especially Van Dyke's matching criterion), even for comparatively simple problems arising from the study of differential equations.

Even in the absence of a formal proof, combining the method of multiple scales and the method of matched asymptotic expansions is clearly an effective approach for some singularly perturbed difference equations, and the work described in this paper expands the range of asymptotic techniques available for difference equations. In the remainder of this section, we outline some important considerations when approaching asymptotic problems in difference equations, and illustrate with further examples a general strategy for seeking asymptotic solutions for difference equations.

%Even in the absence of a formal proof, combining the method of multiple scales and the method of matched asymptotic expansions is clearly an effective approach for some singularly perturbed difference equations, and the work described in this paper expands the range of asymptotic techniques available for difference equations. Moreover, the analysis in this paper enables us to make informed speculations on the features of difference equations that might indicate where the combination of multiple scales and matched asymptotic expansions could be valuable. In the remainder of this section, we outline some important considerations when approaching asymptotic problems in difference equations, and illustrate with further examples a general strategy for seeking asymptotic solutions for difference equations.

In the cases considered in this paper, we began by proposing a multiple scales expansion involving one or two continuum slow time parameters. However, it may not always be necessary to apply the method of multiple scales to a given difference equation, and it is not always obvious how many slow time parameters should be used. Sometimes a regular asymptotic expansion where $x(n)$ is treated as a series in powers of $\eps$ will be sufficient for obtaining an asymptotic solution, as it is for
\begin{equation}
 x(n+1) = \frac{x(n)}{2} + \eps \, x(n)^2, \qquad x(0) = 1.
\end{equation}
In this case, quick inspection reveals that a regular asymptotic expansion will yield each order of $x(n)$ as a sum of terms of the form  $a_j \, 2^{-jn}$ for positive integers $j$, and that secular terms and changes of dominant balance can never arise.

If a regular perturbation expansion fails, it is then important to investigate \emph{why} it fails. In the examples considered in this paper, proposing a regular asymptotic expansion for $x(n)$ at early times would have led to a loss of asymptoticity when $n = \ord(\eps^{-1})$ because of the appearance of secular terms, indicating the need for the method of multiple scales. However, it is also possible to construct examples where a regular asymptotic expansion fails because of a change of dominant balance, indicating a need for matched asymptotic expansions but not the method of multiple scales. Examples of this include the boundary value problems considered in Section 2.8 of \cite{HolmesPert}, or an initial value problem like
\begin{equation}
 x(n+1) = \frac{x(n)}{2} + \frac{\eps}{x(n)}, \qquad x(0) = 1, \label{Nonlinear BL}
\end{equation}
where the rapid decay in the early time solution leads to a change in dominant balance when $x = \ord(\eps^\frac{1}{2})$ and hence $n = \frac{1}{2} \log_2(\frac{1}{\eps}) + \Ord(1)$.

For \eqref{Nonlinear BL}, the fact that the change in dominant balance occurs when $n$ is within $\Ord(1)$ of a logarithmically large critical value indicates that both the early time scaling and the late time scaling will only involve a discrete time parameter; the method of matched asymptotic expansions is needed, while the method of multiple scales is not. In the specific example of \eqref{Nonlinear BL} it is straightforward to perform a full analysis using the method of matched asymptotic expansions because a closed form solution exists for the strongly nonlinear late time equation. In general, however, discrete scale boundary layers in nonlinear problems are likely to lead to strongly nonlinear difference equations for which no closed form solution is known, making a full analysis more difficult.

In cases where the regular perturbation expansion fails because of the appearance of secular terms rather than because of a change of dominant balance, we use the method of multiple scales. If secular terms appear when $n = \ord(\eps^{-1})$, for example, we can introduce a second time scale, $t = \eps \, n$, and construct an asymptotic solution that remains valid while $t = \ord(1)$. However, as $n$ and $t$ grow larger, the two-scale asymptotic expansion may itself fail. One possibility is that the two-scale asymptotic expansion fails because of new secular terms appearing on a longer time scale. This process can be continued indefinitely if each new multiple scales solution fails because of the appearance of secular terms on an even longer time scale.

% to address the secular loss of asymptoticity \emph{before} we reach the stage where we can deal with the change in dominant balance. %Further scales can continue to be introduced until we obtain an asymptotic solution that is uniformly valid for all times, or until there is a failure of asymptoticity associated with a change in dominant balance.

A second possibility for the loss of asymptoticity of a two-scale asymptotic approximation (or indeed any multiple scale asymptotic approximation) is a change in dominant balance. This has been the focus of the work presented in this paper, where we shift from the dominant balance associated with the neighbourhood of a 1-periodic equilibrium (or adiabatic manifold) to the dominant balance associated with the neighbourhood of a 2-periodic equilibrium (or adiabatic manifold). In such cases, we find that the early time multiple scales solution needs to be matched with a late time multiple scales solution; both the method of multiple scales and the method of matched asymptotic expansions are necessary in order to construct a uniformly valid solution.

While we confined our analysis to the discrete logistic equation, it is reasonable to expect similar behaviour in the neighbourhood of any period doubling bifurcation. That is, perturbing a difference equation around a period doubling bifurcation should lead to problems that are amenable to a combined approach based on multiple scales and matched asymptotic expansions. More generally, we expect this behaviour to be characteristic of weakly nonlinear systems of difference equations where the largest eigenvalue of the unperturbed linear problem has unit modulus, but the long time solution involves growth or decay of $x(n)$ leading to a change in dominant balance.

For example, consider the difference equation
\begin{equation}
 x(n+1) - x(n) = - \eps \, x(n) + \frac{\eps^2}{x(n)}, \qquad x(0) = 1, \label{Another Example}
\end{equation}
which is not associated with a period-doubling bifurcation. At early times, the most significant term on the right hand side of \eqref{Another Example} is the $-\eps \, x(n)$ term, which would lead to secular terms if we used a regular perturbation expansion. Using the method of multiple scales to avoid secularity, this leads to exponential decay on the $t = \eps \, n$ time scale. However, the exponential decay in the solution of \eqref{Another Example} cannot go all the way to $x = 0$, because the $\frac{\eps^2}{x(n)}$ term becomes more important as $x$ gets smaller. Hence, there is a change of dominant balance when $x = \ord(\eps^\frac{1}{2})$, or equivalently when $t = \frac{1}{2} \log(\frac{1}{\eps}) + \Ord(1)$, and we can use the method of matched asymptotic expansions to account for the change in dominant balance, precisely as described in this paper.

Moreover, consider the following differential equation, equivalent to Example 3 from Section 1 of \cite{VanHorssen2009}:
\begin{equation}
  x(n+1) = x(n) + \eps \, x(n)^2, \qquad x(0) = \theta > 0.
 \label{NonlinearBlowupDifferenceEqn}
\end{equation} 
In \cite{VanHorssen2009} this is used as an example to show the importance of discreteness, since the difference equation \eqref{NonlinearBlowupDifferenceEqn} has qualitatively different behaviour from the corresponding differential equation: \eqref{NonlinearBlowupDifferenceEqn} has a solution for all $n$, while the corresponding differential equation exhibits finite time blow-up.

However, this can also be interpreted as a problem where there is a change in dominant balance that necessitates the application of matched asymptotic expansions. Applying the method of multiple scales to \eqref{NonlinearBlowupDifferenceEqn} with a continuum slow time variable $t = \eps \, n$ yields
\begin{equation}
 x(n,t) = \frac{\theta}{1 -  \theta \, t} + \eps \, \frac{\theta^2 \, \log(1 - \theta \, t)}{(1 - \theta \, t)^2} + \Ord(\eps^2).
 \label{NonlinearBlowupEarlySol}
\end{equation}
Not only does this solution blow up as $t \rightarrow \theta^{-1}$, we also see that there is a loss of asymptoticity as $t \rightarrow \theta^{-1}$. There is a new dominant balance when $x = \ord(\eps^{-1})$ and hence $t = \theta^{-1} - \ord(\eps)$, enabling us to define late time variables and apply the method of matched asymptotic expansions.

Unlike the the period doubling bifurcations analysed in detail in Sections \ref{S:StaticAlgebraic} and \ref{S:Dynamic} (and, indeed, unlike the example of \eqref{Another Example} described above), the late time scaling for \eqref{NonlinearBlowupDifferenceEqn} leads to a strongly nonlinear difference equation that is not known to have a closed form solution. In practice, it is therefore difficult to make a great deal of progress with the analysis of \eqref{NonlinearBlowupDifferenceEqn}. In principle, however, we see that \eqref{NonlinearBlowupDifferenceEqn} is amenable to the method of multiple scales with a continuum slow time parameter if combined with the method of matched asymptotic expansions in order to obtain a uniformly valid solution. The discrepancy between the finite time blow up of \eqref{NonlinearBlowupEarlySol} and the fact that \eqref{NonlinearBlowupDifferenceEqn} has a solution for all $n$ is not a reason to reject the use of a continuum slow time parameter. Instead, it signals the importance of applying the method of matched asymptotic expansions when there is a failure of asymptoticity associated with a change of dominant balance.

\subsection{Conclusions}

The methods of multiple scales and matched asymptotic expansions are two important asymptotic techniques with a long history of being applied to differential equations. In the present study, we have demonstrated that these two methods can be combined to obtain uniformly valid asymptotic approximations to singularly-perturbed difference equations. By combining these methods, we obtained novel asymptotic approximations to the solutions of autonomous and non-autonomous forms of the logistic map, given in \eqref{BaseEqn-Static} and \eqref{0.0 Dynamic Temp} respectively, when the initial state is close to the unstable period-1 manifold. The resultant asymptotic approximations were presented in \eqref{1.0 composite 2term} for the static problem, and \eqref{Dynamic Composite} for the dynamic problem.

Additionally, the analysis presented in this paper has a number of unusual features that indicate that difference equations may provide very fertile ground for further asymptotic analysis. In each case, for example, we found that the relationship between the early-time independent variables and the late-time independent variables was affine, not linear. That is, converting from $(n,t)$ variables to $(m,s)$ variables involved both shifting and stretching, whereas in most applications of the method of matched asymptotic expansions to differential equations, the equivalent relationship would purely involve stretching. This is because the problems that we considered were not `boundary layer' problems in the classical sense: there was no equivalent to a small parameter multiplying the highest derivative, and the solution did not change rapidly to satisfy a boundary condition. Despite this, the problems analysed did involve a change in dominant balance, and we found that this could successfully be accounted for by finding the new distinguished limit and applying the method of matched asymptotic expansions.

A further significant feature of our analysis is the method used for obtaining the late time scaling in Sections \ref{S:StaticAlgebraicFailure} and \ref{S:DynamicFailure}, in which we considered the doubled map in order to exclude the effects of the fast time scale, which was irrelevant to the new distinguished limit. Without using this method, it would have been very easy to obtain the wrong scaling for the late time problem and thereby fail to recover a uniformly valid solution. We know of no other work where the details of a multiple scales approximation are required in order to find the appropriate distinguished scalings for applying the method of matched asymptotic expansions; seeking similar problems (and equivalent solution methods) in the asymptotic analysis of differential equations is left as an open problem.

\section{Acknowledgements}
The authors would like to acknowledge helpful comments and suggestions from Dr Mohit Dalwadi, who read drafts of this manuscript. We are also grateful for the contributions of the anonymous referees, who provided helpful comments and suggestions. \\ CJL is supported by Australian Laureate Fellowship grant no.~FL120100094 from the Australian Research Council.

\appendix
\section{General form of $X_r(n,t)$ in the case where $\bif = 3+\eps$ and $x(0) = \frac{2}{3}$}
\label{static_appendix}

In the early time analysis of the case where $\bif = 3 + \eps$ and $x(0) = \frac{2}{3}$, it is possible to construct a simple algorithm for obtaining each subsequent term in the asymptotic expansion of $X(n,t)$ in order to satisfy \eqref{1.0 ScaledStatic2}. In order to obtain these terms, we first note that \eqref{1.0 Xn1} gives the solution for $X_0(n,t)$ in the form
\begin{equation}
  X_{0}(n,t) = f_0(t) + g_0(t) \, (-1)^{n},
\end{equation}
with $f_0(t) = \frac{1}{9}$ and $g_0(t) = -\frac{1}{9} \, \e^t$. By continuing the matching process to $\mathcal{O}(\eps^r)$, we may obtain a general form for $X_r$
\begin{align}
 \nonumber\mathcal{O}(\eps^r): \quad
  X_{r}(n+1,t) +&  X_{r}(n,t) = 
  - \frac{1}{3} X_{r-1}(n,t) 
  - 3 \sum_{k=0}^{r-1} X_{k}(n,t) \, X_{r-1-j}(n,t)\\
 & - \sum_{k=0}^{r-2} X_{k}(n,t) \, X_{r-2-k}(n,t) 
  - \sum_{j=1}^{r}\frac{1}{j!} \, \diff{^j X_{r-j}(n+1,t)}{t^j}, \label{1.0 GeneralTerms}
\end{align}

Based on these results, it is possible to show by induction that $X_r(n,t)$ can generally be expressed in the form
\begin{equation}
  X_{r}(n,t) = f_r(t) + g_r(t) \, (-1)^{n}, \label{1.0 Xk Form}
\end{equation}
for any $r$.

To see this, we begin by assuming that $X_R(n,t)$ has the form in \eqref{1.0 Xk Form} for all $R < r$. ,Since we will use the method of multiple scales to assert that the coefficient of the $(-1)^r$ term on the right-hand side of \eqref{1.0 GeneralTerms} is zero for all $r$, we find that \eqref{1.0 GeneralTerms} becomes
\begin{multline}
  X_{r}(n+1,t) + X_{r}(n,t) = 
  - \frac{1}{3} f_{r-1}(t)
  - 3 \sum_{k=0}^{r-1} \Big[ f_{k}(t) \, f_{r-1-k}(t) + g_{k}(t) \, g_{r-1-k}(t) \Big]  \\ 
  - \sum_{k=0}^{r-2}     \Big[ f_k(t)   \, f_{r-2-k}(t) + g_k(t) \, g_{r-2-k}(t) \Big]
  - \sum_{j=1}^{r}\frac{1}{j!} f_{r-j}^{(j)}(t),
\end{multline}
which has a solution of the form
\begin{equation}
 X_r(n,t) = f_{r}(t) + g_r(t) \, (-1)^{n},
\end{equation}
where 
\begin{multline}
 \label{1.0 fk soln}
 f_r(t) = 
  - \frac{1}{6} f_{r-1}(t)
  - \frac{3}{2} \sum_{k=0}^{r-1} \Big[ f_{k}(t) \, f_{r-1-k}(t) + g_{k}(t) \, g_{r-1-k}(t) \Big]  \\ 
  - \frac{1}{2} \sum_{k=0}^{r-2} \Big[ f_{k}(t) \, f_{r-2-k}(t) + g_{k}(t) \, g_{r-2-k}(t) \Big] 
  - \frac{1}{2} \sum_{j=1}^{r}\frac{1}{j!} f_{r-j}^{(j)}(t).
\end{multline}

In order to obtain an expression for $g_r(t)$ in terms of $f_R(t)$ and $g_R(t)$ where $R < r$, and $f_r(t)$, we need to consider the secularity condition at $\mathcal{O}(\epsilon^{r+1})$. Collecting the $(-1)^n$ terms at this order, we find that we require
\begin{equation}
 0 = -\frac{g_r(t)}{3} 
   - 6 \sum_{k=0}^{r}  f_{k}(t) \, g_{r-k}(t)   
   - 2 \sum_{k=0}^{r-1} f_k(t) g_{r-1-k}(t)
   + \sum_{j=1}^{r+1}\frac{1}{j!} g_{r+1-j}^{(j)}(t). 
\end{equation}

Noting that $f_0(t) = \frac{1}{9}$, this leads to a first order differential equation for $g_r(t)$:
\begin{equation}
 \label{1.0 gk ode}
 g_r'(t) -  g_r(t)
  =  6 \sum_{k=1}^{r}  f_{k}(t) \, g_{r-k}(t)   
  + 2 \sum_{k=0}^{r-1} f_j(t) g_{r-1-k}(t)
  - \sum_{j=1}^{r}\frac{1}{j!} g_{r-j}^{(j+1)}(t),
\end{equation}
which must be solved subject to the initial condition $g_r(0) = - f_r(0)$, so that $X_r(0,0) = 0$.

Having developed a systematic way of obtaining the functions $f_r(t)$ and $g_r(t)$, it is straight forward to obtain the asymptotic behaviour of $f_r(t)$ and $g_r(t)$ as $t \to \infty$. Specifically, we find that
\begin{align}
 \label{1.0 fg asymptotic}
 f_r(t) &\sim (-1)^r \, \alpha_r \, \e^{2r t} &
 g_r(t) &\sim (-1)^{r+1} \, \beta_r \, \e^{[2r+1] t},
\end{align}
where $\alpha_r$ and $\beta_r$ are all positive constants.

Again, this can be demonstrated by induction. We begin by noting that \eqref{1.0 fg asymptotic} holds for $r = 0$. Next, we assume that $f_R$ and $g_R$ have this asymptotic form for all $R < r$. Now, considering only the largest terms on the right hand side of  \eqref{1.0 fk soln}, it follows that
\begin{equation}
  f_r(t) \sim
  (-1)^{r} \, \frac{3}{2} \, \sum_{k=0}^{r-1} \Big[ \alpha_{r} \, \alpha_{r-1-k} + \beta_{k} \, \beta_{r-1-k} \Big] \, \e^{2 r t},
\end{equation}
and so
\begin{equation}
 \alpha_r = \frac{3}{2} \, \sum_{k=0}^{r-1} \Big[ \alpha_{k} \, \alpha_{r-1-k} + \beta_{k} \, \beta_{r-1-k} \Big],
\end{equation}
which is clearly positive.

Similarly, assuming the inductive hypothesis and considering only the largest terms on the right hand side of \eqref{1.0 gk ode} leads to
\begin{equation}
  g_r'(t) -  g_r(t)
  \sim  6 \, (-1)^{r+1} \, \sum_{k=1}^{r}  \alpha_r \, \beta_{r-k} \, \e^{[2r+1]t}, 
\end{equation}
and so \eqref{1.0 fg asymptotic} is satisfied with
\begin{equation}
 \beta_r = \frac{6}{2r} \, \sum_{k=1}^{r}  \alpha_r \, \beta_{r-k},
\end{equation}
which is again clearly positive.

More generally, induction can be used to show that 
\begin{align}
 \label{1.0 fg form}
 f_r(t) &= \sum_{k=0}^{r} p_{kr}(t) \, \e^{2r t} &
 g_r(t) &= \sum_{k=0}^{r} q_{kr}(t) \, \e^{[2r +1]t},
\end{align}
where $p_{kr}(t)$ and $q_{kr}(t)$ are both polynomials in $t$.

As $t \to \infty$, the exponentials will determine the dominant behaviour, and we therefore find that
\begin{equation}\label{1.0 asympform}
 X_r(n,\,t) \sim \beta_r (-1)^{r+1+n} \, \e^{[2r+1]t}.
\end{equation}

Recalling that $x(n,\,t) = \sum \epsilon^k \, X_k(n,\,t)$, we note that all of the terms in the series expansion for $x(n,\,t)$ are of the same asymptotic size when $t = \frac{1}{2} \log (\frac{1}{\epsilon}) + \mathcal{O}(1)$ and hence $\epsilon^k \, X_k = \mathcal{O}\big(\epsilon^{-\frac{1}{2}}\big)$. This indicates that the na{\"{\i}}ve series solution fails to be asymptotic when $t = \frac{1}{2} \log (\frac{1}{\epsilon} )+ \mathcal{O}(1)$ and we will need to use an alternative approach.

\bibliographystyle{spmpsci}
\bibliography{DifferenceEqns}

\end{document}